\documentclass{article}

\usepackage{amssymb} 
\usepackage{amsmath}
\usepackage{mathtools} 
\usepackage{color}
\usepackage{tabularx}
\usepackage{scrextend} 
\usepackage{dashrule}
\usepackage{hyperref} 
\usepackage{afterpage} 

\makeatletter
\newcommand\cellwidth{\TX@col@width}
\makeatother

\newcommand{\eref}[1]{\mbox{Eq.\hspace{.12cm}(\ref{eq:#1})}} 
\newcommand{\BigO}[1]{\ensuremath{\operatorname{O}\bigl(#1\bigr)}}

\definecolor{plot_blue}{RGB}{53, 146, 177}
\definecolor{plot_yellow}{RGB}{253, 209, 7}
\definecolor{plot_red}{RGB}{214, 7, 37}

\title{A Class of Exponential Integrators based on Spectral Deferred Correction
    \thanks{This research was supported in part by funding from the Applied Mathematics Department at the University of Washington and NSF grant DMS-1216732.}
}
\author{Tommaso Buvoli 
    \thanks{Department of Applied Mathematics, 
            University of Washington, 
            Lewis Hall, Box 353925, Seattle, WA 98195. 
            (buvoli@uw.edu)}
}
\date{May 28, 2014}

\begin{document}
\maketitle	

\begin{abstract}
We introduce a new class of arbitrary-order exponential time differencing methods based on spectral deferred correction (ETDSDC) and describe a simple procedure for initializing the requisite matrix functions. We compare the stability and accuracy properties of our ETDSDC methods to those of an existing implicit-explicit spectral deferred correction scheme (IMEXSDC). We find that ETDSDC methods have larger accuracy regions and comparable stability regions. We conduct numerical experiments to compare ETD and IMEX spectral deferred correction schemes against a competing fourth-order ETD Runge-Kutta scheme. We find that high-order ETDSDC schemes are the most efficient in terms of function evaluations and overall speed when solving partial differential equations to high accuracy. Our results suggest that high-order ETDSDC schemes are well-suited to work in conjunction with spectral spatial methods or other high-order spatial discritizations. Additionally, ETDSDC schemes appear to be immune to severe order reduction, a problem which affects other ETD and IMEX schemes, including IMEXSDC.
\end{abstract}




\section{Introduction}

In this paper we present a new class of arbitrary-order exponential time differencing (ETD) methods for solving nonlinear evolution equations of the form
\begin{equation*}
\phi_t = \Lambda \phi + \mathcal{N}(t,\phi)
\end{equation*}
where $\Lambda$ is a stiff linear operator and $\mathcal{N}$ is a nonlinear operator. Such systems commonly arise when discretizing nonlinear wave equations including Burgers', nonlinear Schr\"{o}dinger, Korteweg-de Vries, Kuramoto, Navier-Stokes, and the quasigeostrophic equation. ETD Adams methods \cite{beylkin1998ELP,hochbruck2010exponentialreview},  ETD Runge-Kutta methods \cite{cox2002ETDRK4,KassamTrefethen05ETDRK4,krogstad2005IF,hochbruckostermann2005ETDRKSTIFFA,hochbruckostermann2005ETDRKSTIFFB, koikari2005rooted,hochbruck2010exponentialreview}, and ETD general linear methods \cite{ostermann2006general,hochbruck2010exponentialreview} are well-understood, and many of these schemes perform competitively when integrating nonlinear evolution equations \cite{grooms2011IMEXETDCOMP, KassamTrefethen05ETDRK4, loffeld2013comparative}. Despite these advances, no practical high-order exponential integrators have been developed. High-order ETD Adams methods are largely unusable due to their small stability regions, and there are no ETD Runge-Kutta schemes of order greater than five.

Nevertheless, high-order exponential integrators could prove useful if paired with spatial spectral discretizations, especially on periodic domains. Spectral methods exhibit exceptional accuracy and have been shown to be remarkably successful when applied to nonlinear wave equations \cite{fornberg1998PSM,trefethen2000SMM,boyd2013CFSM}. When applying spectral methods on PDEs with smooth solutions, the time integrator often limits the overall order of accuracy.
The development of stable, high-order integrators will allow for more accurate numerical simulations at reduced computational costs and will better balance spatial and temporal accuracy.

In order to develop high-order ETD schemes, we turn our attention to spectral deferred correction methods (SDC), originally developed by Dutt, Greengard, and Rokhlin \cite{Dutt2000SDC}. SDC methods are a class of high-order, self-starting time integrators for solving ordinary differential equations. By pairing Euler's method with a Picard integral equation, SDC methods achieve an arbitrary order of accuracy and favorable stability properties. Remarkably, they are simple to implement, even at high order. In the past decade, there has been a continuing effort to analyze and improve these methods \cite{speck2013multi, tang2013high,  hansen2011order, christlieb2009spectral, huang2006accelerating, liu2008strong, layton2005implications,  guttel2013efficient}. In particular, Minion introduced implicit-explict spectral deferred correction schemes (IMEXSDC) for integrating stiff semilinear systems \cite{Minion2003IMEX}. To date, these methods remain the only practical arbitrary-order IMEX integrators. 
   
   In this paper, we present a new exponential integrator based on spectral deferred correction methods. Our new integrator, which we call ETDSDC, allows for an arbitrary-order of accuracy, has favorable stability properties, and outperforms state-of-the-art ETD schemes when low error tolerances are required. In Section \ref{sec:sdc}, we provide a brief introduction to spectral deferred correction methods before deriving our ETDSDC method and discussing IMEXSDC. In Section \ref{sec:stability_accuracy}, we analyze and compare the stability and accuracy regions of these two methods. In Section \ref{sec:etd_coefficients}, we discuss two techniques for accurately initializing the coefficients for our ETDSDC method. Finally, in Section \ref{sec:numerical_experiments}, we perform numerical experiments comparing our ETDSDC method against IMEXSDC and ETDRK4, a well-known fourth-order exponential integrator \cite{cox2002ETDRK4}.


\section{Spectral Deferred Correction Methods}
\label{sec:sdc}

In this section, we provide a review of Euler-based spectral deferred correction methods  \cite{Dutt2000SDC}, before deriving our ETDSDC method in Section \ref{subsec:etdsdc} and the IMEXSDC method \cite{Minion2003IMEX} in Section \ref{subsec:imexsdc}. To introduce SDC methods, we consider a first-order initial value problem of the form
\begin{equation}
\begin{split}
& \phi'(t) = F(t,\phi) \\
& \phi(a) = \phi_a
\end{split}
\label{eq:model_ode_vanilla}
\end{equation}
where $\phi\in\mathbb{C}^d$ and $F(t,\phi)$ is $\nu$ times differentiable for $\nu \gg 1$. We then shift our attention to a semi-linear first-order initial value problem of the form
\begin{eqnarray}
\begin{split}
&\phi'(t) = \Lambda \phi + \mathcal{N}(t,\phi) \\
& \phi(a) = \phi_a
\end{split}
\label{eq:model_ode_semilinear}
\end{eqnarray}
where again $\phi\in\mathbb{C}^d$, $\mathcal{N}\in C^\nu$, and $\Lambda$ is a $d\times d$ matrix (not necessarily diagonal). The continuity conditions on $\mathcal{N}(t,\phi)$ and $F(t,\phi)$ are stronger than the Lipschitz continuity required for existence and uniqueness, but they ensure that high-order methods can be applied successfully.


\subsection{Preliminaries}
\label{sec:preliminaries}

Spectral deferred correction schemes iteratively improve the accuracy of
an approximate solution to \eref{model_ode_vanilla} by repeatedly solving an integral equation that governs error. This integral equation is of the form
    \begin{equation}
    y(t) = y(a) + \int^t_a g(s,y(s))ds + r(t),
    \label{eq:picard_like_eqn}
    \end{equation}
where $r(a)=0$. As first proposed by Dutt et al. \cite{Dutt2000SDC}, we can approximate the solution to \eref{picard_like_eqn} at points $t_0$, $t_1$, $\ldots$, $t_m$ using the implicit ($\ell=1$) or explicit ($\ell=0$) Euler-like method
    \begin{equation}
    y(t_{n+1}) = y(t_n) + h_n g(t_{n+\ell},y(t_{n+\ell})) + r(t_{n+1}),
    \label{eq:el_method}
\end{equation}
where $h_n=t_{n+1} - t_n$. 

To arrive at the error equation of the form (\ref{eq:picard_like_eqn}), we let $\phi^k(t)$ be an approximate solution to \eref{model_ode_vanilla}, and let the error be ${E(t) = \phi(t) - \phi^k(t)}$.
By considering the integral form of \eref{model_ode_vanilla}, one arrives at
    \begin{equation*}
    \phi(t) = \phi(a) + \int_{a}^{t} F(s,\phi(s)) ds.
    \label{eq:sdc_vanilla_intergral}
    \end{equation*}
Substituting $\phi(t) = \phi^{k}(t) + E(t)$ leads to the integral equation
    \begin{equation}
    E(t) = - \phi^k(t) + \phi^k(a) + E(a) + \int^{t}_{a} F(s,\phi^k(s) + E(s))ds.
    \label{eq:sdc_vanilla_error}
\end{equation}
Introducing the residual
    \begin{equation}
    R(t,a,\phi^k) = \left[ \phi^k(a) + \int^{t}_{a} F(s,\phi^k(s)) ds \right] -\phi^k(t)
    \label{eq:sdc_vanilla_residual_equation}
    \end{equation}
allows us to rewrite \eref{sdc_vanilla_error} as
    \begin{equation}
    E(t) = E(a) + \int^t_a G(s,E(s)) ds + R(t,a,\phi^k),
    \label{eq:sdc_vanilla_correction_equation} \\
    \end{equation}
where
\begin{equation}
G(s,E(s)) = F(s,\phi^k(s) + E(s)) - F(s,\phi^k(s)).
\label{eq:G_definition}
\end{equation}
Rewriting \eref{sdc_vanilla_error} in this manner isolates the residual and the error terms and leads to an equation of the form (\ref{eq:picard_like_eqn}). The residual $R(t,a,\phi^k)$ depends only on known quantities and can be approximated to arbitrary accuracy via numerical quadrature of the function $F(t,\phi^{k}(t))$. If we consider a single timestep of method (\ref{eq:el_method}) applied to \eref{sdc_vanilla_correction_equation}, and suppose that $\phi^k(t)$ is a sufficiently good approximation so that
    \begin{equation*}
    \sup_{t\in[t_{n+1},t_n]}\|E(t)\|  = O(h^m) \hspace{1em} \text{for} \hspace{1em} h=t_{n+1} - t_n \text{ and } m\in\mathbb{N},
    \end{equation*}
then, since $F(t,\phi)$ is Lipchitz continuous in $\phi$, we have that
    \begin{equation*}
    \|h G(s,E(s))\| = h\| F(s,\phi^k(s) + E(s)) - F(s,\phi^{k}(s)) \| = O(h^{m+1}).
    \end{equation*}
Thus, the Euler-like method (\ref{eq:el_method}) is sufficient for estimating $E(t)$ to \mbox{${O(h^{m+1})}$} in the interval $[t_n,t_{n+1}]$. This approximate error, which we denote by $E^k(t)$,  can be used to obtain an $O(h^{n+1})$ accurate solution \mbox{$\phi^{k+1}(t) = \phi^k(t) - E^k(t)$}. This process can be repeated $M$ times to obtain a sequence of increasingly accurate approximations to \eref{model_ode_vanilla}.

To implement this strategy numerically, Dutt et al. proposed to divide each timestep $[t_n,t_{n+1}]$ into $N$ substeps or quadrature nodes which we denote via ${t_{n,1}, \ldots, t_{n,N}}$ \cite{Dutt2000SDC}. This enables us to represent the approximate solution $\phi^k(t)$ as an interpolating polynomial which passes through the quadrature points. 
We can then calculate a provisional solution $\phi^1(t)$ at each node using either forward or backward Euler, and obtain a sequence of higher-order approximations $\phi^{k}(t_{n,j}) = \phi^{k-1}(t_{n,j}) + E^{k-1}(t_{n,j})$ by repeatedly approximating the error $E(t)$ at each quadrature node using (\ref{eq:el_method}).

The choice of the nodes $t_{n,1},\ldots,t_{n,N}$ affects the quality of the quadrature approximation used to determine \eref{sdc_vanilla_residual_equation}. Dutt et al. use Gauss-Legendre points, and Minion has studied the implications of using different quadrature nodes \cite{layton2005implications}. After $M$ correction sweeps, the order of accuracy at each node is $\min(N,M+1)$, regardless of the choice of quadrature nodes \cite{hansen2011order, tang2013high}. 

To simplify our discussion, we consider only a single timestep of spectral deferred correction from $t_n=0$ to $t_{n+1}$. We find it most convenient to describe SDC methods in terms of normalized quadrature points which reduce to the quadrature points if the stepsize $h=1$. Throughout the rest of this paper we will make extensive use of the following definitions:
    \begin{center}
    \renewcommand{\tabcolsep}{10pt}
    \renewcommand{\arraystretch}{1.5}
    \begin{tabular}{llll}
    Stepsize:  & $h = t_{n+1} - t_{n}$ & Normalized nodes: & $\tau_i = t_{n,i}/h$\\
    Substeps: & $h_i = t_{n,i+1} - t_{n,i}$ & Normalized substeps: & $\eta_i = h_i/h$ \\
    \end{tabular}
    \vspace{0.5em}
    \end{center}
We will use the notation SDC$_N^M$ to denote a spectral deferred correction method which uses the quadrature points $\{\tau_i\}_{i=1}^N$, and performs $M$ correction sweeps. For brevity we also use the variables $\phi_i^k = \phi^k(t_{n,i})$ and $E_i^k = E^k(t_{n,i})$ to denote the approximate solution and the error at the $i$th quadrature node after $k$ correction sweeps. 


\subsection{Euler-Based Spectral Deferred Correction Methods} 
\label{subsec:sdc_intro}

We now describe Euler-based spectral deferred correction methods in detail. Implicit and Explicit SDC methods use Implicit or Explicit Euler respectively to determine the provisional solution $\phi^1(t)$ at the quadrature points $h\tau_i$.
Applying the Euler-like method (\ref{eq:el_method}) to \eref{sdc_vanilla_correction_equation} one obtains an approximation of the error $E(t)$ at each of the quadrature points. Every step of this Euler-like method requires approximating the residual term; we describe this process below.

\vspace{1em}
{\bf \noindent Approximating the Residual Term}:  During the $k^{\text{th}}$ correction sweep, $\phi^{k}(t)$ is known at the quadrature points. The residual term (\ref{eq:sdc_vanilla_residual_equation}) can be approximated for $t=h\tau_{i+1}$ and $a=h\tau_i$ at the cost of $N$ function evaluations $F(h\tau_i,\phi_i^k)$ via
    \begin{equation*}
    \hat{R}(h\tau_{i+1},h\tau_{i},\phi^k) = \phi^k(h\tau_i) -\phi^k(h\tau_{i+1}) + I^{i+1}_{i}(\phi^k) 
    \end{equation*} 
where $I_{i}^{i+1}(\phi^k)$ denotes the $N$th order numerical quadrature approximation to
    \begin{equation}
    \int^{h\tau_{i+1}}_{h\tau_{i}} F(s, \phi^{k}(s)) ds.
    \label{eq:sdc_vanilla_quadrature}
    \end{equation}
The coefficients for this numerical quadrature can be obtained for general quadrature points using an algorithm which we propose in Section \ref{sec:etd_coefficients}. For Chebyshev quadrature points, a fast $O(N\log(N))$ matrix-free algorithm exists for computing (\ref{eq:sdc_vanilla_quadrature}) \cite{norris1999spectral}.

\vspace{1em}
\noindent Given the initial condition $\phi_1^1 = \phi(a)$, we can express a single timestep of an SDC$_N^M$ method algorithmically:
 
\begin{center}
\vspace{0.5em}
\renewcommand{\arraystretch}{1.5}
\begin{tabularx}{\textwidth}{|X|}
  \hline
  {\bf Implicit $(\ell = 1)$ or Explicit $(\ell = 0)$ SDC$_\text{N}^\text{M}$ } \hfill Note: $E_1^k$ = 0 \\ \hline

        {\it \textbullet\ Initial Solution (Euler):}
        
        \hspace{1em} \text{{\bf for} i=1 to N-1}
        
        \hspace{2em} $\phi^1_{i+1} = \phi^1_{i} + h_{i} F(h\tau_{i+\ell}, \phi^1_{i+\ell})$
        
        \vspace{1em}
        {\it \textbullet\ Correction \& Update:}
        
        \hspace{1em} \text{{\bf for} k=1 to M}
        
        \hspace{2em} \text{{\bf for} i=1 to N-1}
        
        \hspace{3em} $E_{i+1}^k = E_{i}^k + h_i G(h\tau_{i+\ell}, E^k, \phi^k) + \hat{R}(h\tau_{i+1},h\tau_{i},\phi^k)$

        \noindent \hspace{3em} $\phi^{k+1}_{i+1} = \phi^{k}_{i+1} + E^k_{i+1}$ \vspace{.25em}  \\  \hline
\end{tabularx}
\vspace{0.5em}
\end{center}

\noindent By substituting the expression for $E^k_{i+1}$ into the update formula for $\phi^{k+1}_{i+1}$, noting that $\phi^{k+1}_i = \phi^k_i + E^k_i$, and using \eref{G_definition}, one arrives at the following direct update formula:
    \begin{eqnarray*}
    \phi^{k+1}_{i+1} = \phi^{k+1}_{i} + h_i\left[F(h\tau_{i+\ell},\phi_{i+\ell}^{k+1}) - F(h\tau_{i+\ell}, \phi_{i+\ell}^k) \right] + I_i^{i+1}(\phi^k).
    \end{eqnarray*}
This compact form for spectral deferred correction methods was first mentioned in \cite{Minion2003IMEX} but was not recommended due to potential numerical rounding errors. However, in our numerical experiments, we find that this compact formula leads to simpler codes and equally accurate results. We therefore make use of this compact update formula in all of our codes.


\subsection{ETD Spectral Deferred Correction Methods}
\label{subsec:etdsdc}

We now introduce a new class of exponential integrators based on spectral deferred correction for solving \eref{model_ode_semilinear}, which we repeat here for convenience:
    \begin{equation*}
    \begin{split}
    &\phi'(t) = \Lambda \phi + \mathcal{N}(t,\phi), \\
    & \phi(a) = \phi_a.
    \end{split}
    \end{equation*} 
To derive ETD spectral deferred correction schemes, we seek an error equation of the form
    \begin{equation}
    y(t) = y(a)e^{\Lambda(t-a)} + \int^t_a e^{\Lambda(t-s)} g(s,y(s))ds + r(t).
    \label{eq:etd_picard_like_eqn}
    \end{equation}
We propose to approximate the solution to \eref{etd_picard_like_eqn} by replacing $g(s,y(s))$ with a one-point approximation, leading to the explicit ($\ell=0$) or implicit ($\ell = 1$) ETD Euler-like method
    \begin{equation}
    y(t_{n+1}) = y(t_n) e^{h \Lambda} + \Lambda^{-1} \left[ e^{h \Lambda}-I \right] g(t_{n+\ell},y(t_{n+\ell})) + r(t_{n+1}).
    \label{eq:etd_el_method}
    \end{equation}
To arrive at an error equation of the form (\ref{eq:etd_picard_like_eqn}), we let $\phi^k(t)$ be an approximate solution of \eref{model_ode_semilinear}, and define the error to be $E(t) = \phi(t) - \phi^k(t)$. Applying variation of constants, we obtain the integral form of \eref{model_ode_semilinear},
    \begin{equation*}
    \phi(t) = \phi(a) e^{\Lambda (t-a)} + \int_{a}^{t} e^{\Lambda (t-s)} \mathcal{N}(s,\phi(s)) ds.
    \label{eq:sdc_etd_intergral}
    \end{equation*}
Substituting $\phi(t) = \phi^{k}(t) + E(t)$ leads to the integral equation
    \begin{equation}
    E(t) = -\phi^{k}(t) + \left(\phi^{k}(a) + E(a)\right) e^{\Lambda(t-a)} + \int_{a}^t e^{\Lambda (t-s)}  \mathcal{N}(s, \phi^{k}(s)+E(s)) ds.
    \label{eq:sdc_etd_error}
    \end{equation}
Introducing the residual
    \begin{equation}
    R_{e}(t,a,\phi^k) = \left[ \phi^k(a)e^{ \Lambda (t-a)} + \int_{a}^t e^{\Lambda (t-s)} \mathcal{N}(s, \phi^k(s)) ds \right]  -\phi^k(t)
    \label{eq:sdc_etd_residual_equation}
    \end{equation}
allows us to rewrite \eref{sdc_etd_error} as
    \begin{equation}
    E(t) = E(a)e^{ \Lambda (t-a)} + \int_{a}^t e^{\Lambda (t-s)} H(s,E(s)) ds + R_{e}(t,a,\phi^k),
    \label{eq:sdc_etd_correction_equation}
    \end{equation}
where 
\begin{equation}
H(s,E(s)) = \mathcal{N}(s, \phi^{k}(s)+E(s)) - \mathcal{N}(s, \phi^{k}(s)).
\label{eq:H_definition}
\end{equation}
Now that we have obtained an error equation of the form (\ref{eq:etd_picard_like_eqn}), we are free to proceed in the same manner as Euler-based spectral deferred correction. The provisional solution $\phi^1(t)$ is calculated at the quadrature points using either implicit or explicit ETD Euler and the error at each quadrature point is estimated using (\ref{eq:etd_el_method}). As before, we describe the computation of the residual term.

\vspace{1em}
{\bf \noindent Approximating the Residual Term}:  During the $k^{\text{th}}$ correction sweep, $\phi^{k}(t)$ is known at the quadrature points. The residual (\ref{eq:sdc_etd_residual_equation}) can be approximated for ${t=h\tau_{i+1}}$ and ${a=h\tau_i}$ at the cost of $N$ function evaluations via
    \begin{equation}
    \hat{R}_e(h\tau_{i+1},h\tau_{i},\phi^k(t)) = \phi^{k}(h\tau_i) e^{ h_i \Lambda} -\phi^k(h\tau_{i+1}) + W^{i+1}_{i}(\phi^k)
    \end{equation} 
where $W_{i}^{i+1}(\phi^k)$ denotes the weighted $N$ point numerical quadrature approximation to
    \begin{equation}
    \int_{h\tau_i}^{h\tau_{i+1}} e^{\Lambda (h\tau_{i+1}-s)} \mathcal{N}(s, \phi^{k}(s)) ds
    \label{eq:etd_quad_approx}
    \end{equation}
where the weight function is $w(s) = e^{\Lambda (\tau_{i+1}-s)}$. We describe in detail how to obtain the coefficients for this weighted quadrature in Section \ref{sec:etd_coefficients}. 

\vspace{1em}
\noindent We use ETDSDC$_N^M$ to denote an ETD spectral deferred which performs $M$ correction sweeps on the quadrature points $\{\tau_i\}_{i=1}^N$. Given the initial condition $\phi_1^1=\phi(a)$, we can express a single timestep of an ETDSDC$_N^M$ method algorithmically:
\begin{center}
\vspace{0.5em}
\renewcommand{\arraystretch}{1.5}
\begin{tabularx}{\textwidth}{|X|}
    \hline 
    {\bf Implicit $(\ell = 1)$ or Explicit $(\ell = 0)$ ETDSDC$_\text{N}^\text{M}$} \hfill Note: $E_1^k$ = 0 \\
    \hline

        {\it \textbullet\ Initial Solution (ETD Euler):}
        
        \hspace{1em} \text{{\bf for} i=1 to N-1}
        
        \hspace{2em} $\displaystyle \phi^1_{i+1} = \phi^1_{i} e^{h_i \Lambda} + \Lambda^{-1}\left[ e^{h_i \Lambda}-I \right] \mathcal{N}(h\tau_{i+\ell}, \phi^1_{i+\ell})$
        
        \vspace{1em}
        {\it \textbullet\ Correction \& Update:}
        
        \noindent \hspace{1em} \text{{\bf for} k=1 to M}
        
        \noindent \hspace{2em} \text{{\bf for} i=1 to N-1}
        
        \noindent \hspace{3em} $\displaystyle E_{i+1}^k = E_{i}^k e^{h_i \Lambda} + \Lambda^{-1} \left[ e^{h_i \Lambda}-I \right] H(h\tau_{i+\ell}, E^k, \phi^k) + \hat{R}_e(h\tau_{i+1},h\tau_{i},\phi^k)$
        
        \noindent \hspace{3em} $\phi^{k+1}_{i+1} = \phi^{k}_{i+1} + E^k_{i+1}$
        \vspace{.25em} \\
    \hline 
\end{tabularx}
\vspace{0.5em}
\end{center}

\noindent By substituting the expression for $E^k_{i+1}$ into the update formula for $\phi^{k+1}_{i+1}$, noting that $\phi^{k+1}_i = \phi^k_i + E^k_i$, and using \eref{H_definition}, one arrives at the following direct update formula:
    \begin{align}
    \phi^{k+1}_{i+1} = \phi^{k+1}_{i} e^{h_i \Lambda} + \Lambda^{-1} \left[e^{h_i \Lambda}-1\right] \left[\mathcal{N}(h\tau_{i+\ell},\phi_{i+\ell}^{k+1}) - \mathcal{N}(h\tau_{i+\ell}, \phi_{i+\ell}^k) \right] + W_i^{i+1}(\phi^k).
    \label{eq:etdsdc_compact}
    \end{align}
Though we have derived both an implicit and explicit exponential integrator, we will be solely considering the explicit exponential integrator throughout the rest of this paper.


\subsection{IMEX Spectral Deferred Correction}
\label{subsec:imexsdc}

We now briefly discuss Minion's IMEXSDC$^M_N$ method for solving \eref{model_ode_semilinear} \cite{Minion2003IMEX}. The provisional solution $\phi^1(t)$ is calculated using IMEX Euler. 
The error and residual equations can be derived by repeating the procedure outlined in Section \ref{sec:preliminaries} with ${F(t,y) = \Lambda y + \mathcal{N}(t,y)}$. This leads to 
    \begin{align}
    & E(t) = E(a) + \int_{a}^{t} \left[\Lambda E(s) + G(s,E(s)) \right]ds + R(t,a,\phi^k), \label{eq:sdc_imex_correction_equation} ~~~\\
    & H(s,E(s)) = \mathcal{N}(s,E(s) + \phi^{k}(s)) - \mathcal{N}(s,\phi^{k}(s)) \\
    & R(t,a,\phi^k) = \left[ \phi^{k}(a) + \int_{a}^{t} \left[ \Lambda \phi^{k}(s) + \mathcal{N}(s,\phi^{k}(s)) \right] ds \right]-\phi^{k}(t). \label{eq:sdc_imex_residual_equation}
    \end{align}
Notice that \eref{sdc_imex_correction_equation} is of the form
    \begin{equation}
    y(t) = y(a) + \int^t_a \left[ \Lambda y(s) + g(s,y(s)) \right] ds + r(t).
    \label{eq:imex_picard_like_eqn}
    \end{equation}
We can approximate \eref{imex_picard_like_eqn} by treating the linear term implicitly and the nonlinear term explicitly, yielding the IMEX Euler-like scheme
    \begin{equation*}
    y(t_{n+1}) = (I - h\Lambda)^{-1} \left[y(t_n) + h g(t_{n+\ell},y(t_{n+\ell})) + r(t_{n+1})\right].
    \label{eq:imex_el_method}
    \end{equation*}
The residual term (\ref{eq:sdc_imex_residual_equation}) is approximated exactly as described in Section \ref{subsec:sdc_intro}, except the integrand in  \eref{sdc_vanilla_quadrature} is now ${\Lambda \phi^{k}(s) + N(s,\phi^{k}(s))}$. We denote the quadrature approximation to the residual for IMEXSDC by $\tilde{R}(t,a,\phi)$. Given the initial condition $\phi_1^1=\phi(a)$, we can express a single timestep of an IMEXSDC$_N^M$ method algorithmically:

\begin{center}
\vspace{0.5em}
\renewcommand{\arraystretch}{1.5}
\begin{tabularx}{\textwidth}{|X|}
    \hline 
    {\bf IMEXSDC$_\text{N}^\text{M}$ Method} \hfill Note: $E_1^k$ = 0 \\
    \hline

        {\it \textbullet\ Initial Solution (IMEX Euler):}
        
        \noindent \hspace{1em} \text{{\bf for} i=1 to N-1}
        
        \noindent \hspace{2em} $\displaystyle \phi^1_{i+1} = \left[I-h_i \Lambda \right]^{-1}\left[ \phi^{1}_{i} + h_i N(h\tau_i,\phi^{1}_i)\right]  $
        
        \vspace{1em}
        {\it \textbullet\ Correction \& Update:}
        
        \noindent \hspace{1em} \text{{\bf for} k=1 to M}
        
        \noindent \hspace{2em} \text{{\bf for} i=1 to N-1}
        
        \noindent \hspace{3em} $\displaystyle E_{i+1}^k = \left[I-h_i \Lambda \right]^{-1} \left[E^k_i + h_i H(h\tau_i,E, \phi^{k}) + \tilde{R}(h\tau_{i+1},h\tau_i,\phi^k) \right]$
        
        \noindent \hspace{3em} $\phi^{k+1}_{i+1} = \phi^{k}_{i+1} + E^k_{i+1}$
        \vspace{.25em} \\
    \hline  
\end{tabularx}
\vspace{0.5em}
\end{center}

\noindent By rewriting the error formula implicitly so that
    \begin{equation*}
    E_{i+1}^k =  \left[E^k_i + h_i(\Lambda E^k_{i+1} + H(h\tau_i,E, \phi^{k})) + \tilde{R}(h\tau_{i+1},h\tau_i,\phi^k) \right],
    \end{equation*}
substituting this expression into the update formula for $\phi^{k+1}_{i+1}$, and noting that
    \begin{equation*}
     E_{i+1}^k = \phi^{k+1}_{i+1} - \phi^{k}_{i+1}, \hspace{2em} \phi^{k+1}_i = \phi^k_i + E^k_i
    \end{equation*}
one arrives at the following direct update formula:
    \begin{equation*}
    \phi^{k+1}_{i+1} = \left[I-h_i \Lambda \right]^{-1} \left[\phi_{i}^{k+1} - (h_i \Lambda) \phi_{i+1}^{k} + h_i(\mathcal{N}(h\tau_i,\phi^{k+1}_i) - \mathcal{N}(h\tau_i,\phi_{i}^{k}) ) + \tilde{I}_i^{i+1}(\phi^k)\right]
    \end{equation*}
where $\tilde{I}^{i+1}_i(\phi^k)$ denotes the numerical quadrature approximation to
    \begin{equation*}
    \int^{h\tau_{i+1}}_{h\tau_{i}} \Lambda \phi^k(s) + \mathcal{N}(s, \phi^{k}(s)) ds.
    \end{equation*}


\section{Stability and Accuracy}
\label{sec:stability_accuracy} 

Determining the stability properties of IMEX and ETD integrators is non-trivial. A commonly used approach is to consider the model problem
    \begin{equation}
    \begin{split}
    & \phi' = \mu \phi + \lambda \phi \\
    & \phi(0) = 1
    \end{split}
    \label{eq:stability_model_problem}
    \end{equation}
where $\mu,\lambda \in \mathbb{C}$ and the terms $\mu \phi$, $\lambda \phi$ act as the linear and nonlinear term respectively. This model problem highlights stability for \eref{model_ode_semilinear} when it is possible to simultaneously diagonalize both the linear and nonlinear operators around a fixed point. Though this analysis does not extend to general linear systems, it has proven useful for predicting stability properties of IMEX and ETD methods on a variety of partial differential equations \cite{grooms2011IMEXETDCOMP}.  

Applying an ETDSDC$_{N}^{M}$ or IMEXSDC$_{N}^{M}$ method on \eref{stability_model_problem} leads to a recursion relation of the form 
\begin{equation*}
\phi(t_{n+1}) = \psi^M_N(r,z) \phi(t_n)
\end{equation*}
where $r=\mu h$, $z=\lambda h$, and $h$ denotes the timestep. As with all one-step methods, the stability region is defined as
\begin{equation*}
\mathcal{S} = \{(r,z)\in \mathbb{C}^2, |\psi^M_N(r,z)|\le 1 \}.
\end{equation*}
We list the stability functions $\psi^M_N(r,z)$ for ETDSDC$^M_N$ and IMEXSDC$^M_N$ schemes in Table \ref{tab:stab_funs}.

\begin{table}
\begin{center}
    \renewcommand{\arraystretch}{1.5}
    \begin{tabularx}{\linewidth}{X}
    {\bf ETDSDC Stability Functions \hfill $\psi_1^k(r,z) = 1$} \\ \hline
        \vspace{-2em}
        \begin{equation}
        \renewcommand{\arraystretch}{3}
        \begin{array}{l c l}
        \psi_{i+1}^{1} &=& \displaystyle e^{r \eta_i} \psi^{1}_{i} + \frac{e^{r \eta_i} - 1}{r} z \psi_{i}^{1} \\
        \psi_{i+1}^{k+1} &=& \displaystyle e^{r \eta_i} \psi_{i}^{k+1} + \frac{e^{r \eta_i} - 1}{r} z (\psi_{i}^{k+1} - \psi_{i}^{k}) + z\sum_{j=1}^N \mathbf{W}_{i,j} \psi_{j}^{k} \hspace{1em} \\
    \text{where} & & \displaystyle \mathbf{W}_{ij} = \int^{\tau_{i+1}}_{\tau_i} e^{r(\tau_{i+1}-s)} L_j(s) ds, \hspace{2em}
        L_j(s) = \prod^{N}_{\substack{l=1\\l \neq j}} \frac{(s-\tau_l)}{(\tau_j-\tau_l)}.        
        \end{array}
        \label{eq:etdsdc_stability_function}
        \end{equation} \\
    \end{tabularx}
\end{center}
\begin{center}
    \renewcommand{\arraystretch}{1.5}
    \begin{tabularx}{\linewidth}{X}
    {\bf IMEXSDC Stability Functions \hfill $\psi_1^k(r,z) = 1$} \\ \hline
        \vspace{-2em}
        \begin{equation}
        \renewcommand{\arraystretch}{3}
        \begin{array}{l c l}
        \psi^1_{i+1} &=& \displaystyle \left(\frac{1 + z \eta_i}{1-r\eta_i} \right) \psi_{i}^{1} \\
        \psi^{k+1}_{i+1} &=& \displaystyle \left(\frac{\psi_{i}^{k+1} + \eta_iz \left( \psi^{k+1}_{i} - \psi^{k}_{i} \right) - r\eta_i \psi_{i+1}^{k} + (r + z) \sum_{j=1}^N \mathbf{I}_{i,j} \psi_{j}^{k}}{1-r\eta_i} \right) \\
        \text{where} & & \displaystyle  \mathbf{I}_{ij} = \int^{\tau_{i+1}}_{\tau_i} L_j(s) ds, \hspace{2em}
        L_j(s) = \prod^{N}_{\substack{l=1\\l \neq j}} \frac{(s-\tau_l)}{(\tau_j-\tau_l)}.
        \end{array}
        \label{eq:imexsdc_stability_function}
        \end{equation}
    \end{tabularx}
\end{center}
\caption{Stability functions for ETDSDC$^M_N$ and IMEXSDC$^M_N$ methods. As $r \to 0$ the stability functions of both methods limit to that of an explicit SDC$^M_N$ method.}
\label{tab:stab_funs}
\end{table}

We choose to analyze stability for PDEs with linear dispersion and dissipation; thus, $r = h\mu$ and $z = h\lambda$ are complex-valued. Several strategies have been proposed for effectively visualizing the resulting four-dimensional stability region. As in \cite{beylkin1998ELP, cox2002ETDRK4, krogstad2005IF}, we choose to overlay two-dimensional slices of the stability regions, each corresponding to a fixed $r$ value. For the sake of brevity, we focus our attention on 8th order methods where ${N=8},~{M=7}$ and on 16th order methods where ${N=16},~{M=15}$. For all methods, we select the Chebyshev quadrature nodes
\begin{equation*}
\tau_i = \frac{1}{2}\left(1 - \cos\left(\frac{\pi(i-1)}{N-1}\right)\right)\hspace{2em} i=1,\ldots,N.
\end{equation*}

\begin{figure}[htp]
\centering

\par{Stability Region Plots}

\hrulefill \newline
\begin{tabular}{c}
{\tiny \textcolor{plot_blue}{\hdashrule[0.2ex]{2em}{2pt}{}} $r = 0$ \hspace{1em}}
{\tiny \textcolor{plot_yellow}{\hdashrule[0.2ex]{2em}{2pt}{}} $r = R_0/2$ \hspace{1em}}
{\tiny \textcolor{plot_red}{\hdashrule[0.2ex]{2em}{2pt}{}} $r = R_0$} \hspace{1em}
{\tiny \textcolor{black}{\hdashrule[0.2ex]{2em}{2pt}{}} $r = 2R_0$} \\
\end{tabular}

\hrulefill \vspace{1em}

\begin{minipage}[b]{0.24\linewidth}
\centering
\par{\hspace{1em} \tiny ETDSDC$_8^7$}
\end{minipage}
\begin{minipage}[b]{0.24\linewidth}
\centering
\par{\hspace{1em} \tiny IMEXSDC$_{8}^{7}$}
\end{minipage}
\begin{minipage}[b]{0.24\linewidth}
\centering
\par{\hspace{1em} \tiny ETDSDC$_{16}^{15}$}
\end{minipage}
\begin{minipage}[b]{0.24\linewidth}
\centering
\par{\hspace{1em} \tiny IMEXSDC$_{16}^{15}$}
\end{minipage}
\vspace{1em}

\par{\tiny \bf Dissipative Model Problem: $r \in -1 \cdot [0, 30]$ and $R_0 = -30$ \vspace{1em} }

\begin{minipage}[b]{0.24\linewidth}
\centering
\includegraphics[width=1\linewidth]{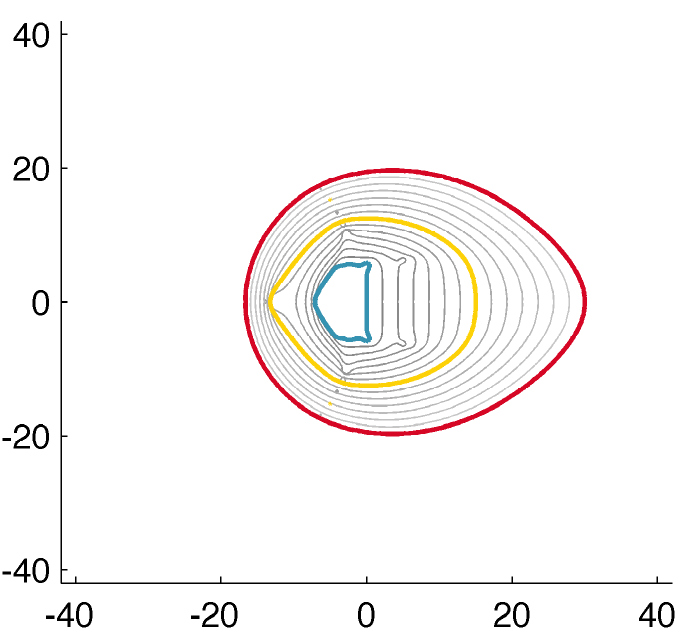}
\end{minipage}
\begin{minipage}[b]{0.24\linewidth}
\centering
\includegraphics[width=1\linewidth]{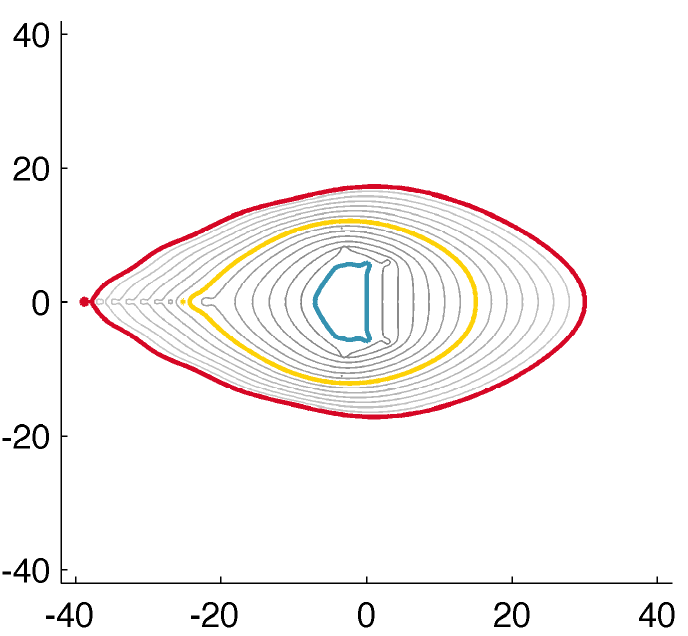}
\end{minipage}
\begin{minipage}[b]{0.24\linewidth}
\centering
\includegraphics[width=1\linewidth]{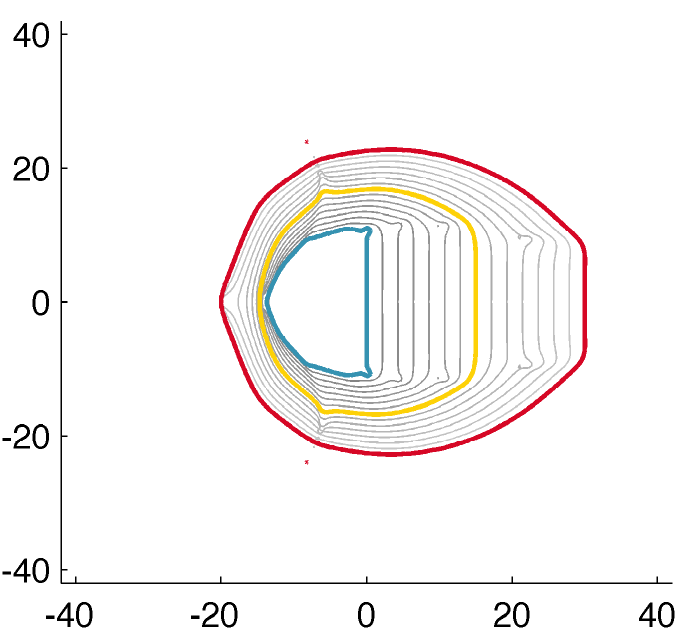}
\end{minipage}
\begin{minipage}[b]{0.24\linewidth}
\centering
\includegraphics[width=1\linewidth]{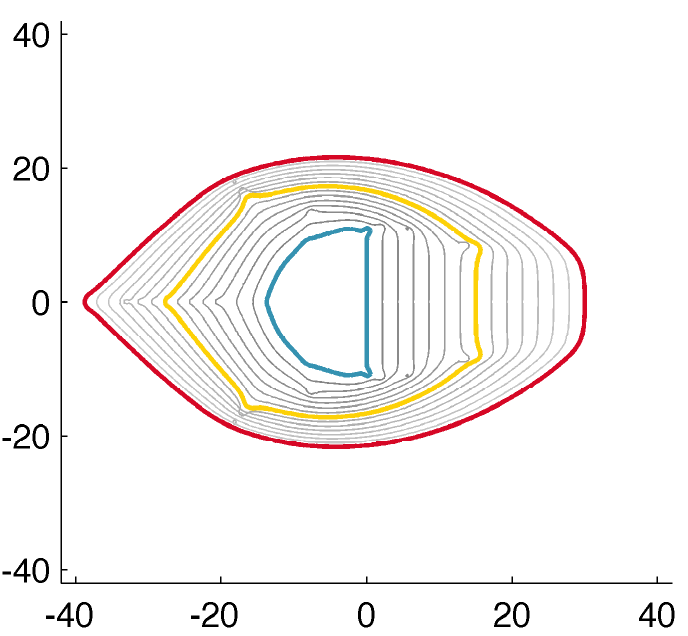}
\end{minipage}

\par{\tiny \bf Dispersive Model Problem: $r \in 1i \cdot [0, 30]$  and $R_0 = 30i$ \vspace{1em} }

\begin{minipage}[b]{0.24\linewidth}
\centering
\includegraphics[width=1\linewidth]{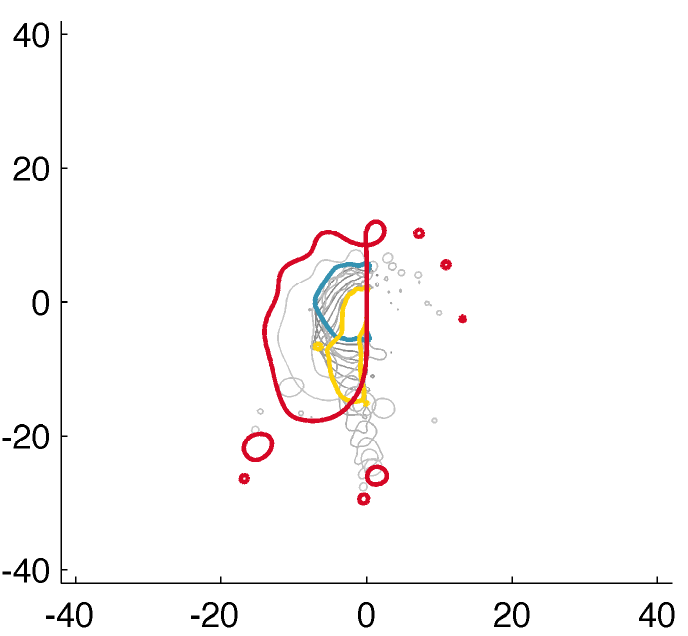}
\end{minipage}
\begin{minipage}[b]{0.24\linewidth}
\centering
\includegraphics[width=1\linewidth]{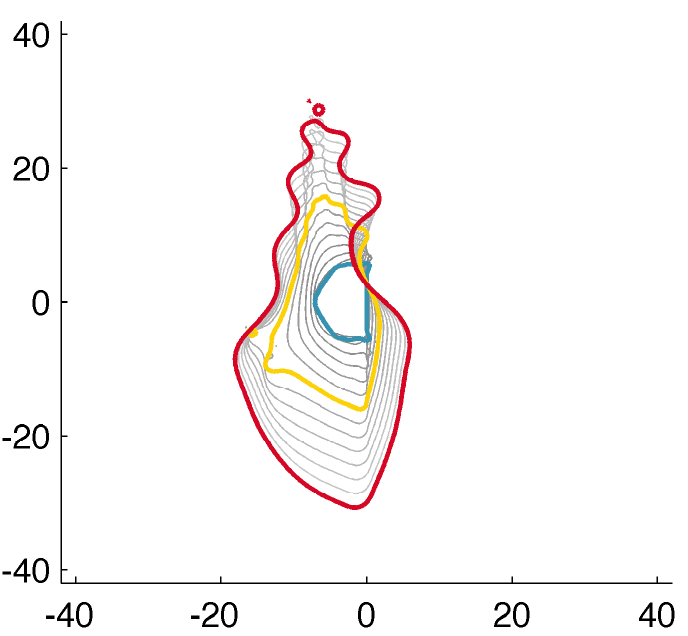}
\end{minipage}
\begin{minipage}[b]{0.24\linewidth}
\centering
\includegraphics[width=1\linewidth]{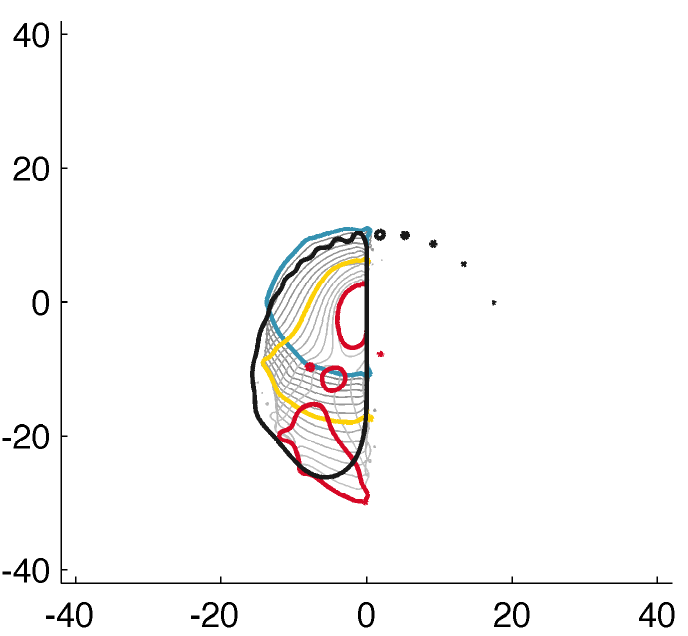}
\end{minipage}
\begin{minipage}[b]{0.24\linewidth}
\centering
\includegraphics[width=1\linewidth]{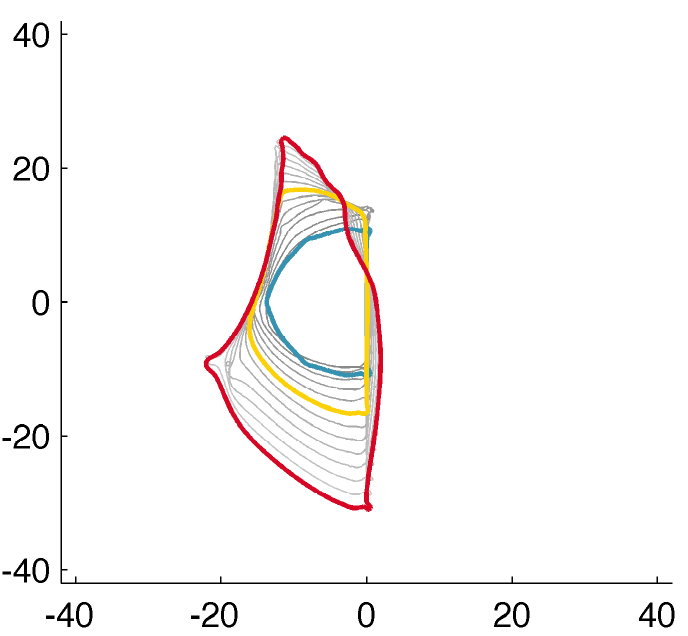}
\end{minipage}

\par{\tiny \bf Dissipative/Dispersive Model Problem: $r \in \exp(3\pi i/4) \cdot [0, 30]$ and $R_0 = 30e^{3\pi i/4}$ \vspace{1em} }

\begin{minipage}[b]{0.24\linewidth}
\centering
\includegraphics[width=1\linewidth]{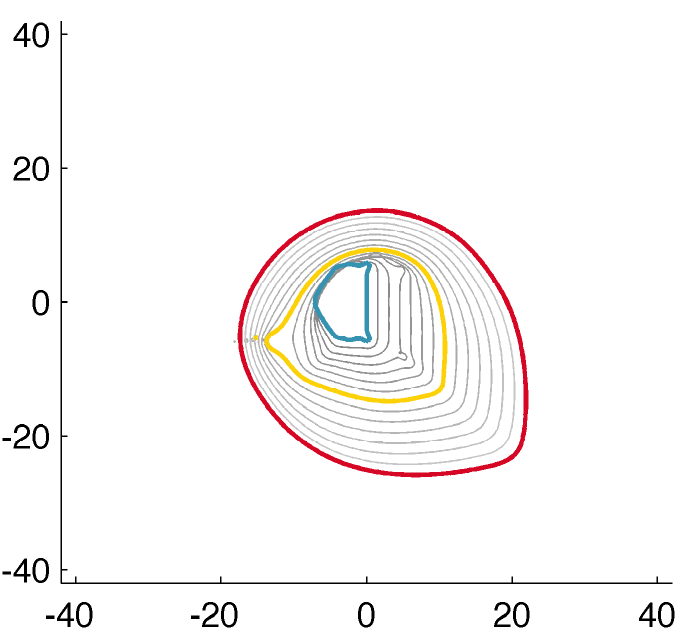}
\end{minipage}
\begin{minipage}[b]{0.24\linewidth}
\centering
\includegraphics[width=1\linewidth]{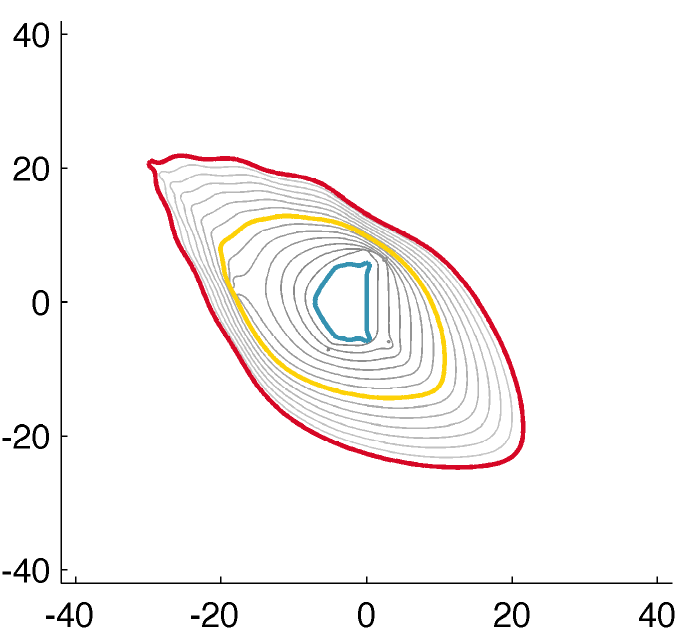}
\end{minipage}
\begin{minipage}[b]{0.24\linewidth}
\centering
\includegraphics[width=1\linewidth]{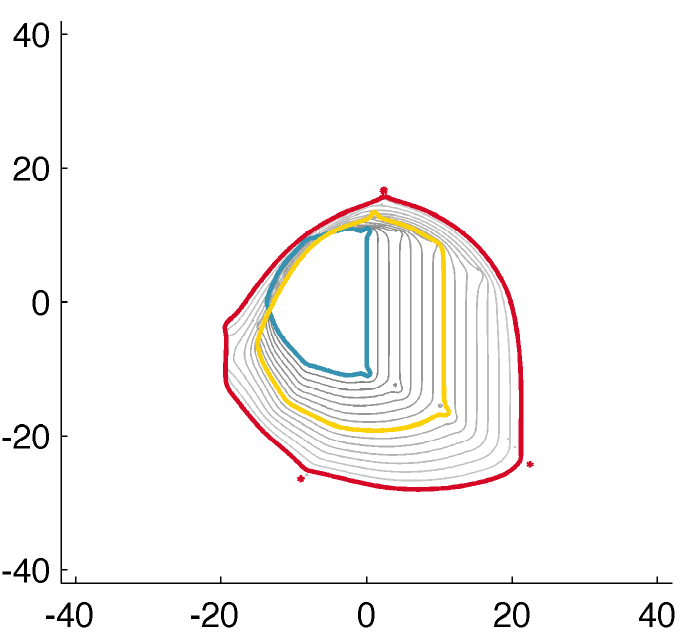}
\end{minipage}
\begin{minipage}[b]{0.24\linewidth}
\centering
\includegraphics[width=1\linewidth]{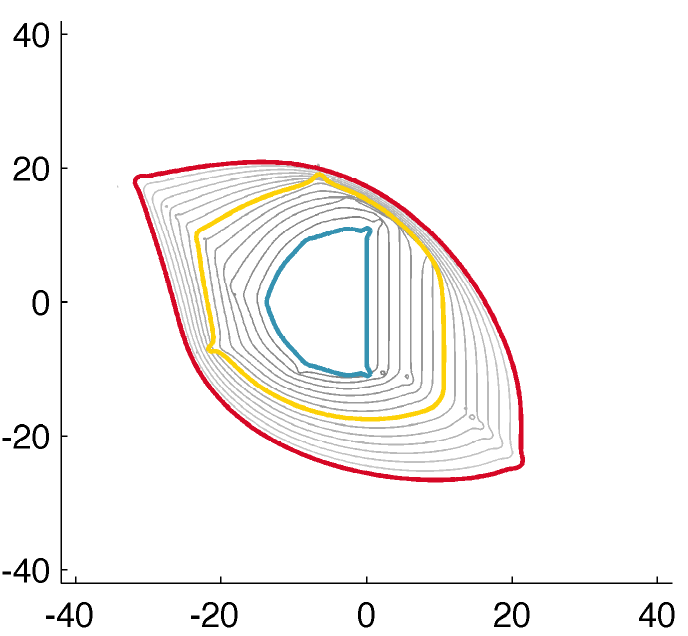}
\end{minipage}

\caption{Stability regions for 8th order and 16th order methods with Chebyshev quadrature nodes. Colored contours correspond to different $r$ values as described in the legend. We plot an additional black contour for the ETDSDC$^{15}_{16}$ method on the dispersive model problem to show that stability regions eventually grow for sufficiently large imaginary $r$. For large $|r|$, increasing the order of the ETD and IMEX methods does not lead to significantly larger stability regions.}
\label{fig:stability}
\end{figure}

We pick a range of real, imaginary, and complex $r$ values to simulate nonlinear PDEs with varying degrees of linear dispersion and dissipation. We plot stability regions pertaining to
\begin{equation}
r \in -1 \cdot [0,30], \hspace{2em} r \in 1i \cdot  [0,30], \text{ and} \hspace{2em} r \in \exp(3\pi i/4) \cdot [0,30]
\label{eq:r_value_range}
\end{equation}
in Figure \ref{fig:stability}. For these three $r$ ranges, we find that the stability regions of all methods grow as $|r|$ increases. For imaginary $r$, the stability regions for ETDSDC methods temporarily decrease before growing. Though all methods exhibit satisfactory stability properties, IMEXSDC methods allow for coarser timesteps on a wider range of $(r,z)$. Overall, our results suggest that both IMEXSDC methods and ETDSDC methods exhibit good stability properties on a wide range of stiff nonlinear evolution equations. 

When analyzing spectral deferred correction methods, it is also common to plot accuracy regions.  Accuracy regions highlight the restrictions on the stepsize $h$ so that error after one timestep is smaller than $\epsilon>0$.  They are simply  defined as
\begin{equation*}
\mathcal{A}_\epsilon = \{(r,z)\in \mathbb{C}^2, |\psi^M_N(r,z) - \exp(r+z)|\le \epsilon \}.
\end{equation*}
They were introduced in \cite{Dutt2000SDC} for comparing the efficiency of high-order methods, and provide a more detailed picture than stability regions which solely differentiate between convergent and divergent $(r,z)$ pairs. 

We find that as $|r|$ increases, the accuracy region containing $z=0$ decreases rapidly for ETDSDC$^M_N$ methods and vanishes entirely for IMEXSDC$^M_N$ methods. This behavior can be understood from \eref{etdsdc_stability_function} and \eref{imexsdc_stability_function}. For the ETDSDC$^M_N$ methods it follows that ${\psi_N^M(r,0) = \exp(r h)}$; moreover, since the stability function $\psi_N^M(r,z)$ is continuous, then for any $\epsilon > 0$, there exists a nontrivial accuracy region surrounding $z=0$. The same cannot be said for IMEXSDC schemes since \eref{imexsdc_stability_function} satisfies the weaker relation ${\psi(r,0) =  \exp(rh) + O(rh)}$; hence, as $r$ becomes sufficiently large, there need not exist an accuracy region around $z=0$.

\begin{figure}[htp]
\centering

\par{Accuracy Region Plots}

\hrulefill \newline
\begin{tabular}{c}
{\tiny \textcolor{plot_blue}{\hdashrule[0.2ex]{2em}{2pt}{}} $r = 0$ \hspace{1em}}
{\tiny \textcolor{plot_yellow}{\hdashrule[0.2ex]{2em}{2pt}{}} $r = R_0/2$ \hspace{1em}}
{\tiny \textcolor{plot_red}{\hdashrule[0.2ex]{2em}{2pt}{}} $r = R_0$} \\
\end{tabular}

\hrulefill

\begin{minipage}[b]{0.24\linewidth}
\centering
\par{\hspace{1em} \tiny ETDSDC$_8^7$}
\end{minipage}
\begin{minipage}[b]{0.24\linewidth}
\centering
\par{\hspace{1em} \tiny IMEXSDC$_{8}^{7}$}
\end{minipage}
\begin{minipage}[b]{0.24\linewidth}
\centering
\par{\hspace{1em} \tiny ETDSDC$_{16}^{15}$}
\end{minipage}
\begin{minipage}[b]{0.24\linewidth}
\centering
\par{\hspace{1em} \tiny IMEXSDC$_{16}^{15}$}
\end{minipage}
\vspace{1em}

\par{\tiny \bf Dissipative Model Problem: $r \in -1 \cdot [0, 30]$ \vspace{1em} }

\begin{minipage}[b]{0.24\linewidth}
\centering
\par{\tiny $R_0 = -5$}
\includegraphics[width=1\linewidth]{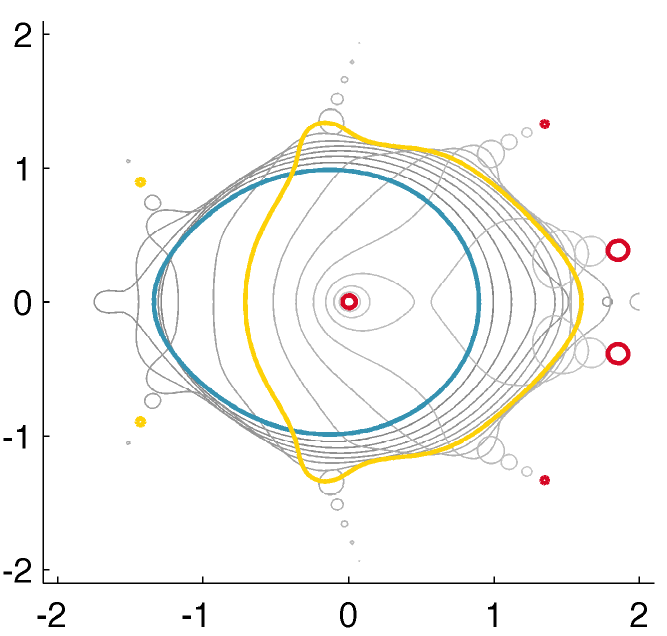}
\end{minipage}
\begin{minipage}[b]{0.24\linewidth}
\centering
\par{\tiny $R_0 = -2$}
\includegraphics[width=1\linewidth]{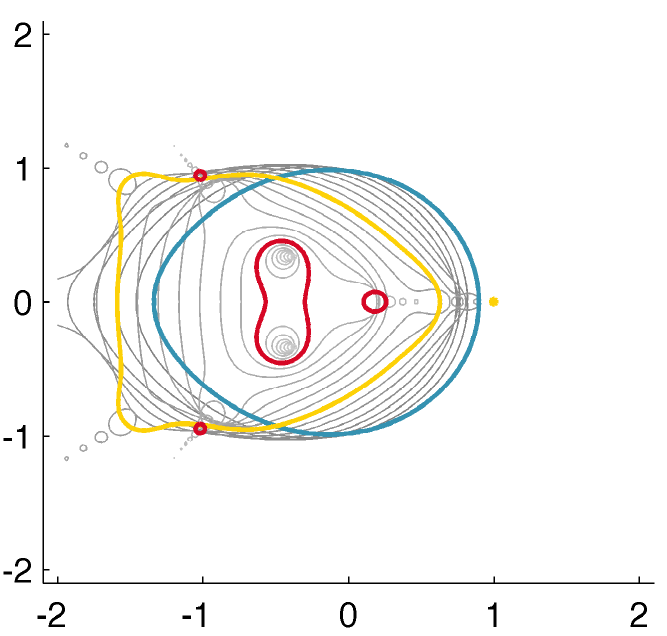}
\end{minipage}
\begin{minipage}[b]{0.24\linewidth}
\centering
\par{\tiny $R_0 = -30$}
\includegraphics[width=1\linewidth]{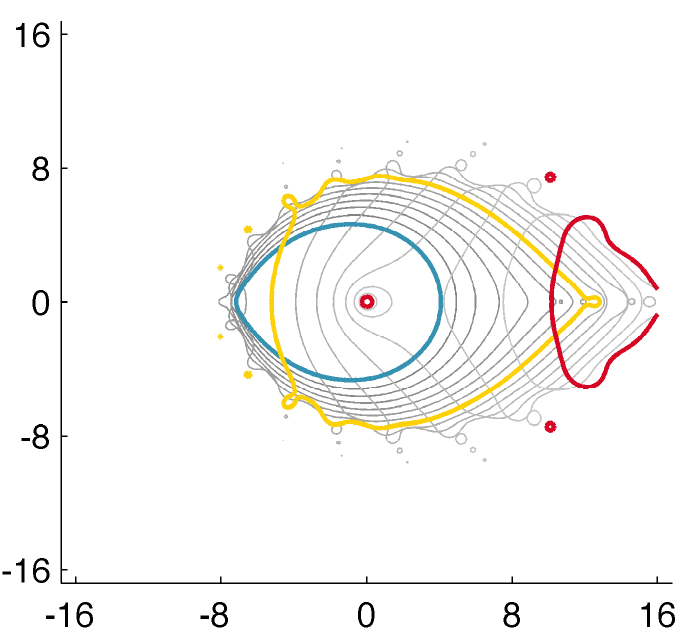}
\end{minipage}
\begin{minipage}[b]{0.24\linewidth}
\centering
\par{\tiny $R_0 = -20$}
\includegraphics[width=1\linewidth]{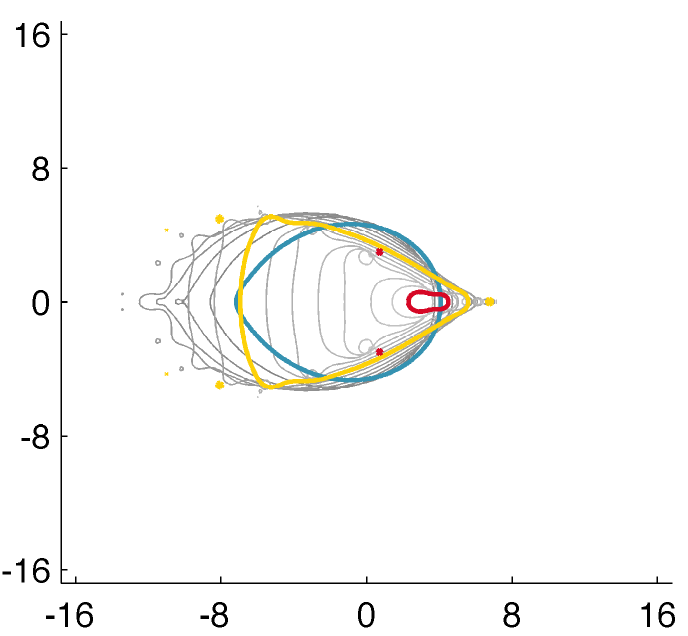}
\end{minipage}

\par{\tiny \bf Dispersive Model Problem: $r \in 1i \cdot [0, 30]$ \vspace{1em} }

\begin{minipage}[b]{0.24\linewidth}
\centering
\par{\tiny $R_0 = 5i$}
\includegraphics[width=1\linewidth]{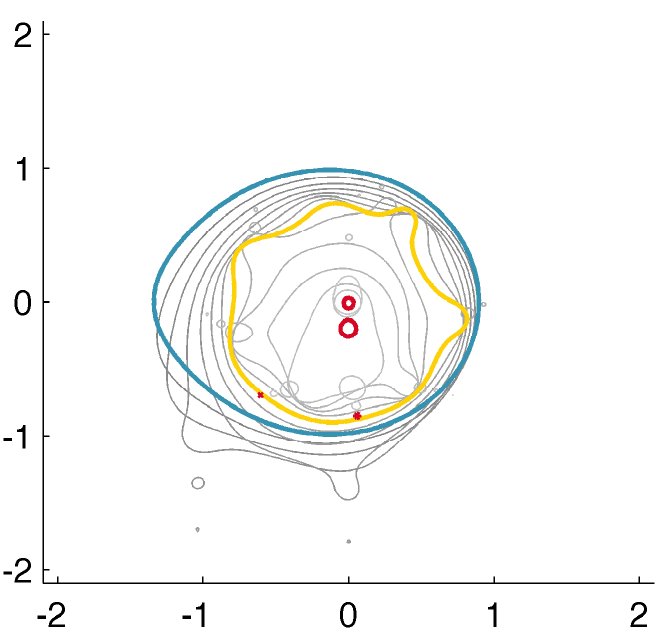}
\end{minipage}
\begin{minipage}[b]{0.24\linewidth}
\centering
\par{\tiny $R_0 = 2i$}
\includegraphics[width=1\linewidth]{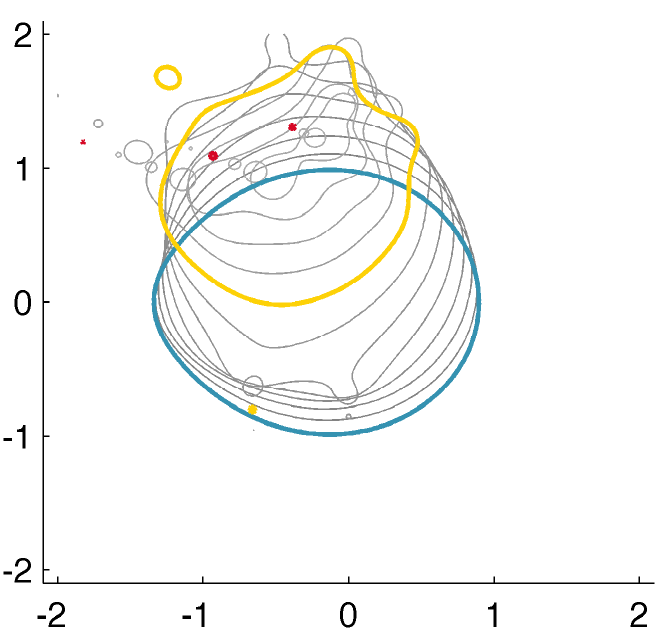}
\end{minipage}
\begin{minipage}[b]{0.24\linewidth}
\centering
\par{\tiny $R_0 = 15i$}
\includegraphics[width=1\linewidth]{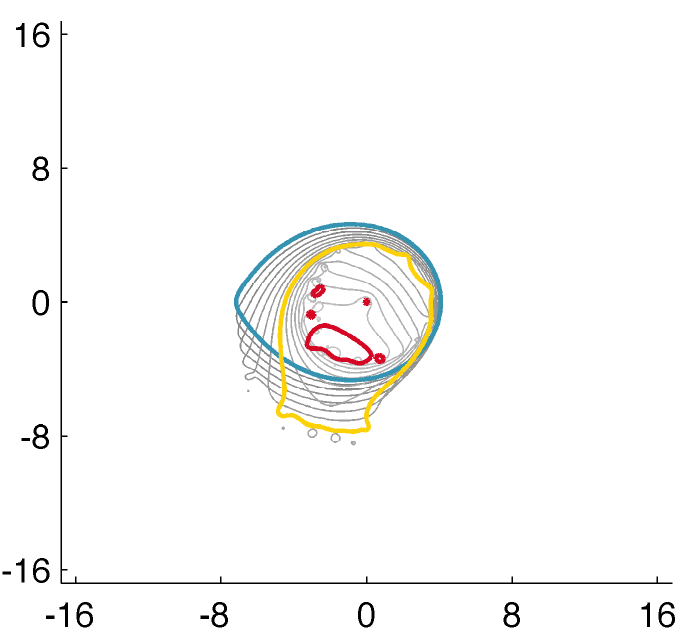}
\end{minipage}
\begin{minipage}[b]{0.24\linewidth}
\centering
\par{\tiny $R_0 = 9i$}
\includegraphics[width=1\linewidth]{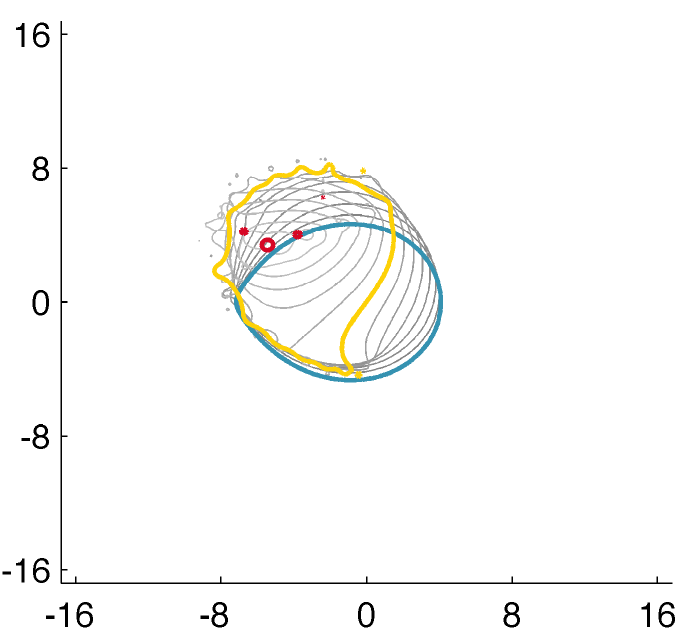}
\end{minipage}

\par{\tiny \bf Dissipative/Dispersive Model Problem: $r \in \exp(3\pi i/4) \cdot [0, 30]$ \vspace{1em} }

\begin{minipage}[b]{0.24\linewidth}
\centering
\par{\tiny $R_0 = 5e^{3i\pi /4}$}
\includegraphics[width=1\linewidth]{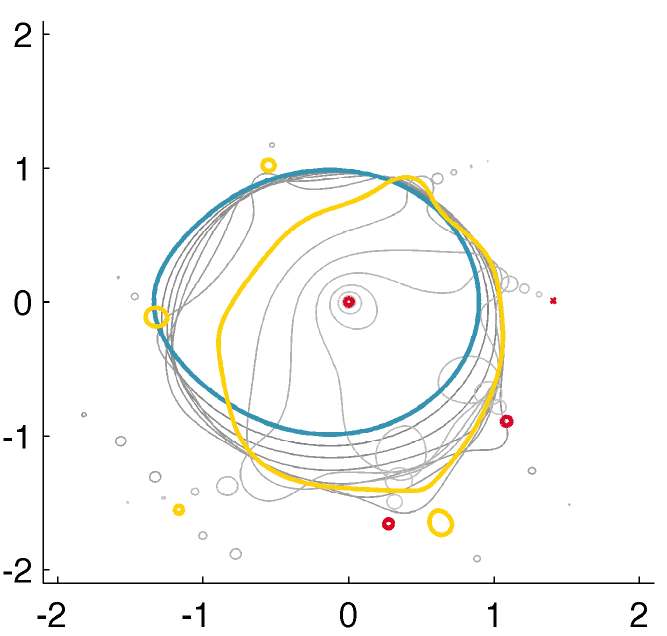}
\end{minipage}
\begin{minipage}[b]{0.24\linewidth}
\centering
\par{\tiny $R_0 = 2e^{3i\pi /4}$}
\includegraphics[width=1\linewidth]{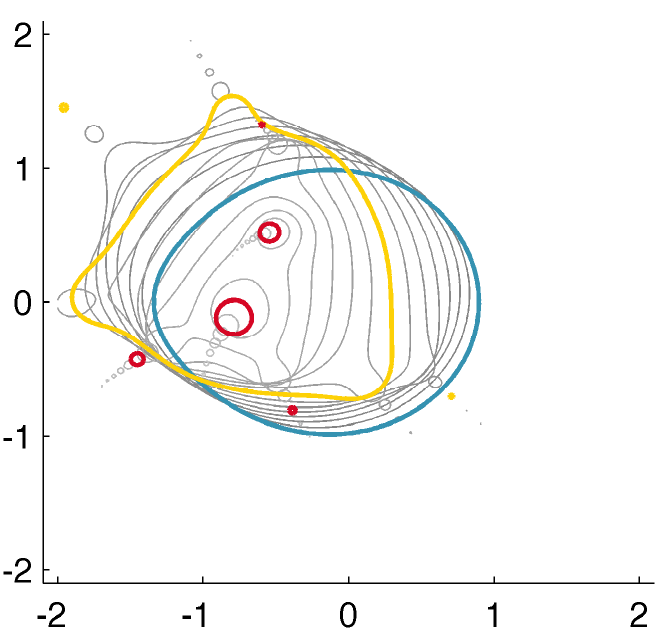}
\end{minipage}
\begin{minipage}[b]{0.24\linewidth}
\centering
\par{\tiny $R_0 = 26e^{3i\pi /4}$}
\includegraphics[width=1\linewidth]{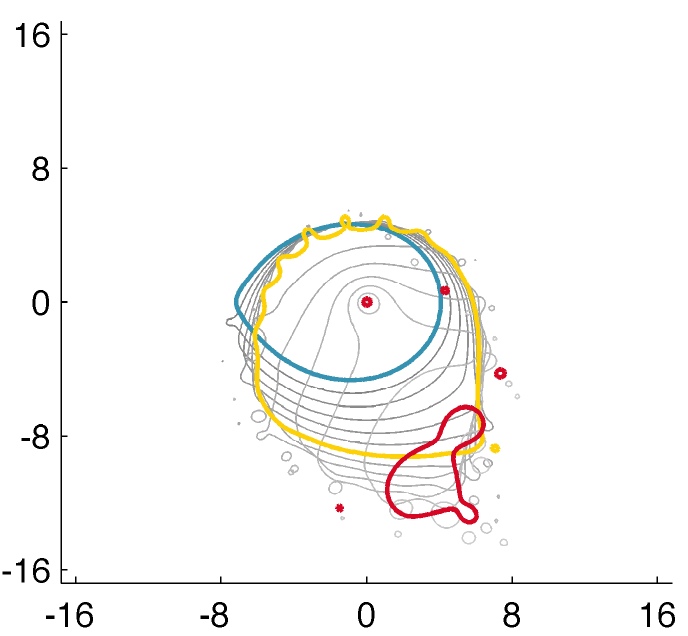}
\end{minipage}
\begin{minipage}[b]{0.24\linewidth}
\centering
\par{\tiny $R_0 = 16e^{3i\pi /4}$}
\includegraphics[width=1\linewidth]{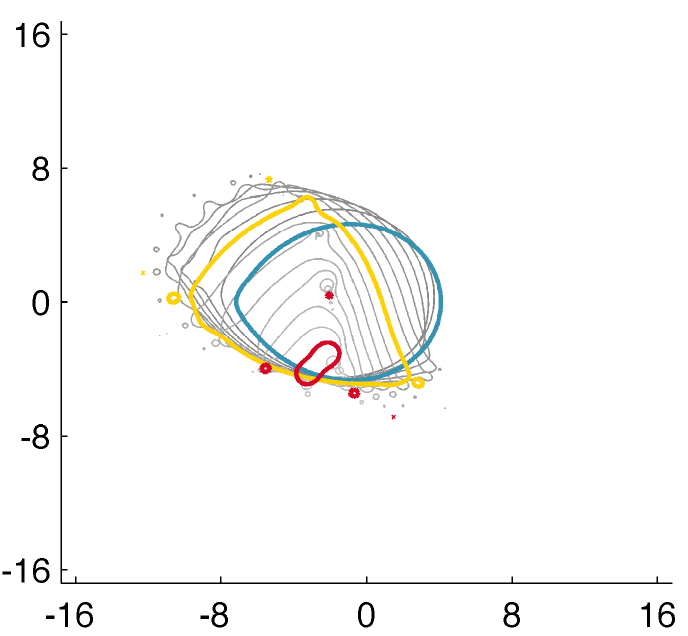}
\end{minipage}

\caption{Accuracy regions corresponding to $\epsilon = 1\times 10^{-8}$ for 8th order and 16th order methods with Chebyshev quadrature nodes. Colored contours correspond to different $r$ values as described in legend. We choose $R_0$ in each figure so that the \textcolor{plot_red}{red} contour marks a near vanishing accuracy region around $z=0$. As expected, 16th order methods possess larger accuracy regions for a wider range of $r$ than 8th order methods.}
\label{fig:accuracy}
\end{figure}

We present accuracy regions for ${\epsilon = 1\times10^{-8}}$ in Figure \ref{fig:accuracy}. We consider the three ranges of $r$ values in (\ref{eq:r_value_range}), but due to rapidly shrinking accuracy regions, we are only able to visualize different subsets of $r$ values for each numerical method. ETDSDC$^M_N$ schemes outperform IMEXSDC$^M_N$ schemes for all tested values. Accuracy regions for the ETD methods decrease more slowly, and the non-vanishing accuracy regions around $z=0$ guarantee accuracy for any $r$ so long as $z$ is chosen sufficiently small. The MATLAB code used to generate these figures can be found in \cite{BuvoliETDZenodo} and can be easily modified to generate stability and accuracy plots for other ETDSDC or IMEXSDC methods.

\section{Calculating $W^{i+1}_{i}(\phi^k)$}
\label{sec:etd_coefficients}
Every iteration of an ETDSDC$^M_N$ method requires computing $W^{i+1}_{i}(\phi^k)$, which denotes the weighted quadrature approximation to 
    \begin{equation*}
    \int_{h\tau_i}^{h\tau_{i+1}} e^{\Lambda (h\tau_{i+1}-s)} \mathcal{N}(s, \phi^{k}(s)) ds.
    \end{equation*}
To arrive at a formula for $W_{i}^{i+1}(\phi^k)$, we let $\mathbf{N}_l(\phi) = \mathcal{N}(h\tau_l, \phi(h\tau_l))$ and replace $\mathcal{N}(s,\phi^k(s))$ in \eref{etd_quad_approx} with the Lagrange interpolating polynomial $L(s)$ that passes through the quadrature points $\{(h\tau_{l}, \mathbf{N}_l(\phi^k)) \}_{l=1}^N$ so that
    \begin{equation}
    W^{i+1}_i(\phi^k)
    = \int^{h\tau_{i+1}}_{h\tau_{i}} e^{\Lambda(h\tau_{i+1} - s)} L(s) ds
    = \sum_{l=1}^N w_{i,l} \hspace{.1em} \mathbf{N}_l(\phi^k). \label{eq:W_lagrange}
    \end{equation} 
For low-order methods, explicit formulae for $w_{i,l}$ can be derived by forming $L(s)$ and repeatedly applying integration by parts. Unfortunately, this direct calculation leads to increasingly involved formulae for large $N$. We therefore seek a general procedure for determining  $w_{i,l}$ for any $N$. We propose to express the weights $w_{i,l}$ in terms of the well-known functions
    \begin{equation*}
    \varphi_n(z) = \frac{1}{(n-1)!} \int^1_0 e^{z(1-\sigma)}\sigma^{n-1} d\sigma.
    \end{equation*}
    using a stable algorithm developed by Fornberg for determining finite difference coefficients \cite{fornberg1988}. We describe our algorithm in Section \ref{sec:w_algorithm}, before discussing $\varphi$ functions and two well-known methods for initializing them in Section \ref{sec:phif}.


\subsection{Proposed Algorithm}
\label{sec:w_algorithm}

To arrive at a convenient expression for $W^{i+1}_{i}(\phi^k)$, we propose to apply the change of variables
    \begin{equation}
    s = h\left[(\tau_{i+1} - \tau_i) \sigma + \tau_i \right] = h_i \sigma + h \tau_i, 
    \label{eq:w_cov}
    \end{equation}
to the integral term in (\ref{eq:W_lagrange}), expand the Lagrange interpolating polynomial $L(s(\sigma))$ as a Taylor polynomial, and rewrite the result in terms of $\varphi$ functions. 
Applying the change of variables (\ref{eq:w_cov}) leads to
    \begin{equation*}
    h_i \int^{1}_{0} e^{h_i \Lambda (1-\sigma)} L(s(\sigma)) d\sigma = h_i \int^{1}_{0} e^{h_i \Lambda (1-\sigma)} P_i(\sigma) d\sigma, 
    \end{equation*}
where $P_i(\sigma)$ is the Lagrange interpolating polynomial which passes through the points
    \begin{equation*}
    \{(q_{i,l}, \mathbf{N}_l(\phi^k)) \}_{l=1}^N 
    \hspace{1em} \text{and} \hspace{1em}
    q_{i,l} = (\tau_l - \tau_i)/(\tau_{i+1} - \tau_{i})
    \end{equation*}
denote the scaled, translated quadrature nodes $h\tau_i$ under the transformation ($\ref{eq:w_cov}$). Next, we define the finite difference coefficients $a^{(i)}_{j,l}$ so that
    \begin{equation*}
    \left. \frac{d^j}{d\sigma^j} P_i(\sigma) \right|_{\sigma=0} = \sum_{l=1}^N a^{(i)}_{j,l} \hspace{.05em} \mathbf{N}_l(\phi^k).
    \end{equation*}
Expanding $P_i(\sigma)$ as a Taylor polynomial we obtain
    \begin{equation*}
    W^{i+1}_{i}(\phi^k) = h_i \int^{1}_{0} e^{h_i \lambda (1-\sigma)} \sum_{j=0}^{N-1} \left[ \frac{\sigma^j}{j!} \sum_{l=1}^{N} a^{(i)}_{j,l} \hspace{0.1em} \mathbf{N}_l(\phi^k) \right] d\sigma.
    \end{equation*}
Reordering terms we arrive at
    \begin{eqnarray*}
    W^{i+1}_{i}(\phi^k) &=& h_i \sum_{l=1}^N \left[ \mathbf{N}_l(\phi^k)  \sum_{j=0}^{N-1} \left[ \frac{a^{(i)}_{j,l}}{j!} \int^{1}_{0}  e^{h_i \Lambda (1-\sigma)} \sigma^j d\sigma \right] \right] \\
    &=& h_i \sum_{l=1}^N \left[ \mathbf{N}_l(\phi^k)  \sum_{j=0}^{N-1} \left[ a^{(i)}_{j,l} \varphi_{j+1}(h_i \Lambda) \right] \right].
    \end{eqnarray*}
By defining the functions
    \begin{equation}
    w_{i,l}(z) = h_i \sum_{j=0}^{N-1} a^{(i)}_{j,l} \varphi_{j+1}(z),
    \label{eq:w_function}
    \end{equation}
we obtain a convenient expression for the weighted quadrature rule:
    \begin{equation*}
    W^{i+1}_{i}(\phi^k) = \sum_{l=1}^N w_{i,l}(h_i\Lambda) \mathbf{N_k}(\phi^k).
    \end{equation*}
To successfully implement this procedure, we must determine the finite difference coefficients $a_{j,l}^{(i)}$ and the matrix functions $\varphi_n(h_i\Lambda)$. The coefficients $a^{(i)}_{j,l}$ can be rapidly obtained using the stable algorithm presented in \cite{fornberg1988}. We define the functions:
\begin{itemize}
    \item $weights(z_0,[q_1,\ldots,q_n],m)$: returns a finite difference matrix $\mathbf{a}$ for computing $m$ derivatives at $z_0$, assuming $q_j$ are the quadrature points. This calling sequence is consistent with the implementation in \cite{fornberg1998classroom}.
    \item $initPhi(z,n)$: returns the functions $\varphi_i(z)$ for $i=0,\ldots,n$. We discuss two possible implementations in Section \ref{sec:phif}.
\end{itemize}
The algorithm for computing $w_{i,l}(z)$ for an ETDSDC$_N^M$ method can be written as: 

\begin{center}
\vspace{0.5em}
\renewcommand{\arraystretch}{1.5}
\begin{tabularx}{\textwidth}{|X|}
    \hline 
    {\bf Computing $w_{i,l}(h_i\Lambda)$} \\ \hline
    \text{{\bf for} i=1 to N}

    \hspace{1em} $[\varphi_0(h_i \Lambda), \ldots, \varphi_N(h_i \Lambda)] = \text{initPhi}(h_i\Lambda,N)$
        
    \hspace{1em}  \text{{\bf for} j=1 to N}     
    
    \hspace{2em} $q_j = (\tau_j - \tau_i)/(\tau_{i+1} - \tau_{i})$ 
    
    \hspace{1em} $a^{(i)} = \text{weights}(0,[q_1, \ldots q_N],N-1)$
    
    \hspace{1em} \text{{\bf for} l=1 to N}
    
    \hspace{2em} \text{{\bf for} j=0 to N-1}
    
    \hspace{3em} $w_{i,l}(h_i\Lambda) = w_{i,l}(h_i\Lambda) + a^{(i)}_{j,l} \varphi_{j+1}(h_i \Lambda)$ \vspace{.5em} \\ \hline  
\end{tabularx}
\vspace{0.5em}
\end{center}

\noindent When computing $w_{i,l}(h_i\Lambda)$, it is convenient to save $\varphi_{0}(h_i\Lambda)$ and $\varphi_{1}(h_i\Lambda)$ since both are required for the ETD Euler method.


\subsection{$\varphi$ Functions}
\label{sec:phif}

The coefficients of all exponential integrators can be expressed in terms of $\varphi$ functions \cite{hochbruck2010exponentialreview,berland2007expint,minchev2005,koikari2007error}. The $n$th $\varphi$ function can be defined in the following ways:
    \begin{align}
    & \text{{\it Integral Form}:} &  \varphi_n(z) & = \begin{dcases} e^{z} & n = 0 \\ \frac{1}{(n-1)!} \int^1_0 e^{z(1-s)} s^{n-1} ds & n> 0 \end{dcases} \\
    & \text{{\it Series Form:}} & \varphi_n(z) & = \sum_{k=0}^\infty \frac{z^{k}}{(k+n)!} \label{eq:phi_taylor_def} \\
    & \text{{\it Recursion Relation:}} & \varphi_n(z) & = \frac{\varphi_{n-1}(z) - \frac{1}{(n-1)!}}{z}, \hspace{1em} \varphi_0(z) = e^z 
    \label{eq:phi_recursion}
    \end{align}
The first few $\varphi_n(z)$ are given by
    \begin{equation*}
    \varphi_0(z) = e^z, \hspace{1em} \varphi_1(z) = \frac{e^z - 1}{z}, \hspace{1em} \varphi_2(z) = \frac{e^z - 1 - z}{z^2}, \hspace{1em} \varphi_3(z) = \frac{e^z - 1 - z - \tfrac{1}{2}z^2}{z^3}.
    \end{equation*}
We can now rewrite the compact update formula (\ref{eq:etdsdc_compact}) as 
    \begin{equation*}
    \phi^{k+1}_{i+1} = \varphi_0(h_i \Lambda) \phi^{k+1}_{i}  + \varphi_1 (h_i \Lambda) \left[N(h\tau_{i+\ell},\phi_{i+\ell}^{k+1}) - N(h\tau_{i+\ell}, \phi_{i+\ell}^k) \right] + W_i^{i+1}(\phi^k).
    \end{equation*}
From their series definition, it follows that the functions $\varphi_n(z)$ are entire; nevertheless, it is well-known that explicit formula for $\varphi_n(z)$ are prone to catastrophic numerical roundoff error for small $|z|$. Various strategies for overcoming this difficulty have been compared extensively  \cite{ashi2009comparison}. We briefly outline a method  based on scaling and squaring \cite{koikari2007error} and a method based on contour integration \cite{KassamTrefethen05ETDRK4}. Other approaches involve Krylov subspace approximations \cite{hochbruck1998exponential,hochbruck1997krylov} and improved contour integrals \cite{trefethen2007} but we do not consider them in this paper.

\subsubsection{Taylor/Pad\'{e} Scaling and Squaring Algorithm}
The scaling and squaring algorithm for calculating $\varphi$ functions is a generalization of a well-known algorithm for computing matrix exponentials \cite{higham2009scaling}. For small $|z|$, $\varphi_n(z)$ can be accurately evaluated via the Taylor series (\ref{eq:phi_taylor_def}) or via the diagonal $(m,m)$ Pad\'{e} approximation, whose explicit formula is given in \cite{skaflestad2009scaling}. This initial approximation can be used to obtain $\varphi_n(z)$ for large $|z|$ by repeatedly applying the well-known scaling relation 
    \begin{equation}
    \varphi_n(z) = \frac{1}{2^n} \left[ \varphi_0 \left(\tfrac{z}{2} \right) \varphi_n \left(\tfrac{z}{2} \right) + \sum_{i=1}^{n} \frac{\varphi_{i}\left(\tfrac{z}{2} \right)}{(n-i)!} \right].
    \label{eq:phi_scaling_relation}
    \end{equation}
We present pseudocode for an $m$-term Taylor series procedure for initializing $\phi_i(\Lambda)$ in Table \ref{tab:phi_methods}. A MATLAB implementation of the Pad\'{e} scaling and squaring algorithm is freely available in \cite{berland2007expint} and can be easily used to initialize $\varphi_n(\Lambda)$ for both scalar and matrix $\Lambda$.

\subsubsection{Contour Integration Algorithm}

An alternative algorithm for initializing ETD coefficients was first suggested in \cite{KassamTrefethen05ETDRK4}. Since the functions $\varphi_n(z)$ are entire, Cauchy's integral formula can be used to obtain $\varphi(z)$ at problematic regions near $z = 0$. We highlight this procedure for both scalar and matrix $\Lambda$ in Table \ref{tab:phi_methods} assuming that the explicit formula for $\varphi_n(z)$ is known. If this is not the case, then it is convenient to combine \eref{phi_recursion} with the discretized contour integral so that
\begin{equation}
\varphi_n(\Lambda) = \frac{1}{P} \sum_{j=0}^{P-1} \frac{\varphi_{n-1}(\Lambda + re^{i\theta}) - 1/(n-1)!}{\Lambda + Re^{i\theta}}.
\label{eq:phi_contour_recursion}
\end{equation}
This allows one to progressively evaluate $\varphi_n(\Lambda)$ for $n=1,\ldots,N$. For scalar $\Lambda$ we use \eref{phi_contour_recursion} when $|\Lambda| < 1$ and \eref{phi_recursion} when $|\Lambda|\ge 1$. For matrix $\Lambda$ we find that the technique based on scaling and squaring is faster and more accurate, especially for matrices with large norm.

\begin{table}
\begin{center}
\renewcommand{\arraystretch}{1.5} 
\begin{tabularx}{\textwidth}{|X|} \hline
    {\bf $m$-Term Taylor Scaling \& Squaring Procedure for Matrix/Scalar $\Lambda$} \\ \hline

    {\it \textbullet\  Select Scaling Factor:}
        
        Let $s\in\mathbb{N}$ so that $\|\Lambda/2^s\|_\infty < \delta(m)$ \hfill
        See Appendix \ref{ap:A1} or \cite{koikari2007error} for choosing $\delta(m)$.
        \vspace{1em}
        
        {\it  \textbullet\  Initialize $\varphi_i(\Lambda/2^s)$ via Horner's Method:} \vspace{0.5em}
        
        {\bf for} i=0 to N  \vspace{.25em}
        
            \hspace{1em} $ P_i = \frac{\Lambda}{(m+i)!} + \frac{\mathbf{I}}{(m+i-1)!}$ \vspace{.5em}
        
            \hspace{1em} {\bf for} k=0 to m-2   \vspace{.25em}
        
                \hspace{2em} $ P_i = \Lambda P_i + \frac{\mathbf{I}}{(m+i-2-k)!}$   \vspace{.25em}
        
                \hspace{1em} $\varphi_i(\Lambda/2^s) = P_i$ \vspace{0.25em}
                
        \vspace{1em}
        {\it  \textbullet\ Obtain $\varphi_i(\Lambda)$ via \eref{phi_scaling_relation}:}
        \vspace{1em}
        
        {\bf for} i=1 to s
    
        \hspace{1em} {\bf for} n=0 to N
    
                \hspace{2em} $\displaystyle \varphi_n(\Lambda/2^{s-i}) = \frac{1}{2^n} \left[ \varphi_0 \left(\Lambda/2^{s-i+1} \right) \varphi_n \left(\Lambda/2^{s-i+1} \right) + \sum_{i=1}^{n} \frac{\varphi_{i}\left(\Lambda/2^{s-i+1} \right)}{(n-i)!} \right]$ \vspace{1em}
                
            \\ \hline

\end{tabularx}

\vspace{-1pt}

    \renewcommand{\arraystretch}{1.5} 
    \renewcommand{\tabcolsep}{0pt}
    \begin{tabularx}{\textwidth}{|X|} 
    {
        \renewcommand{\tabcolsep}{5pt}
        \begin{tabularx}{\cellwidth}{X|X}
            {\bf Contour Integral for Scalar} $|\Lambda| < 1$ & {\bf Contour Integral for Matrix} $\Lambda$ \\ \hline
            { 
                {\it \textbullet\ Cauchy Integral Formula:}
                \vspace{1em}
                
                \begin{addmargin}[1em]{0em}
                    $\displaystyle{\varphi_{n}(\Lambda) = \frac{1}{2\pi i} \oint_\Gamma \frac{\varphi_{n}(z)}{(z - \Lambda)}dz}\vspace{1em}$
                \end{addmargin}
                
                {\it \textbullet\ Choosing $\Gamma$:} \vspace{1em}
        
            \begin{addmargin}[1em]{1em}
            Let $\Gamma = R e^{i\theta}  + \Lambda$ for $\theta\in[0,2\pi]$. The radius $R$ should be chosen so that contour never comes near the origin. \vspace{1em}
            
            $\displaystyle \varphi_n(\Lambda) = \frac{1}{2\pi} \int^{2\pi}_0 \varphi_n(Re^{i\theta} + \Lambda) d\theta$
            \vspace{1em}
            
            \end{addmargin}           
        
            {\it \textbullet\ Discretization via Trapezoidal Rule:} \vspace{1em}
        
            \begin{addmargin}[1em]{0em}
            Let $\theta_j = 2\pi j/P$, then for $P$ sufficiently large, $\varphi_{n}(\Lambda)$ is approximately \newline
            \end{addmargin}

            \hspace{.7em} $\displaystyle{\frac{1}{P} \sum_{j=0}^{P-1} \varphi_{n}(\Lambda + Re^{i\theta_j})}.$
            \vspace{1em}
            \begin{addmargin}[1em]{0em}
            For scalar $|\Lambda|\ge 1$, use \eref{phi_recursion}.
            \end{addmargin}
            } & {
                {\it \textbullet\ Cauchy Integral Formula: \vspace{1em}}
        
                \begin{addmargin}[1em]{0em}
                $\displaystyle{\varphi_{n}(\Lambda) = \frac{1}{2\pi i} \oint_\Gamma \varphi_{n}(z)(z \mathbf{I} - \Lambda)^{-1} dz}$
                \vspace{1em}
                \end{addmargin}
                
                {\it \textbullet\ Choosing $\Gamma$:} \vspace{1em}
                
                \begin{addmargin}[1em]{1em}
                Let $\Gamma = R e^{i\theta} + z_0$ for $\theta\in[0,2\pi]$. The radius $R$ and center $z_0$ must be chosen so that contour encloses spectrum of $\Lambda$. 
                
                \vspace{1em}
                $\displaystyle \varphi_n(\Lambda) = \frac{1}{2\pi} \int^{2\pi}_0 \varphi_n(Re^{i\theta} + \Lambda) d\theta$
                \vspace{1em}
                
                \end{addmargin}
                
                 {\it \textbullet\ Discretization via Trapezoidal Rule:} \vspace{1em}
                 
                \begin{addmargin}[1em]{0em}
                Let $\theta_j = 2\pi j/P$, $\gamma_j = R \exp(\theta_j) + z_0$, then for $P$ sufficiently large, $\varphi_{n}(\Lambda)$ is approximately
                \end{addmargin}
                \vspace{1em}
                        
                \hspace{.7em} $\displaystyle{ \frac{1}{P} \sum_{j=0}^{P-1} \varphi_{n}(\gamma_j) \left(\mathbf{I} + \frac{(z_0 \mathbf{I} - \Lambda)}{Re^{i\theta_j}}\right)^{-1}}\vspace{.5em}$
                
                \hspace{.7em} where $\varphi_n(\gamma_j)$ initialized like scalar $\Lambda$.  \vspace{-1em}
            }
        \end{tabularx}    
    } \\ \hline
    \end{tabularx}
\end{center}

\vspace{1em}
\caption{Scaling \& squaring, and contour integral methodology for initializing $\varphi_n(\Lambda)$.}
\label{tab:phi_methods}
\end{table}


\section{Numerical Experiments}
\label{sec:numerical_experiments}

In this section, we numerically solve four partial differential equations in order to compare ETDSDC$^M_N$ and IMEXSDC$^M_N$ methods of orders $4,8,16,32$ against the fourth-order Runge-Kutta method (ETDRK4) developed in \cite{cox2002ETDRK4}.  We have chosen to include ETDRK4 in our tests since it was shown to perform competitively \cite{grooms2011IMEXETDCOMP, KassamTrefethen05ETDRK4}, and provides a good reference for comparing SDC based schemes to existing ETD and IMEX methods. We provide our MATLAB and Fortran implementation of ETDSDC$^M_N$, IMEXSDC$^M_N$ and ETDRK4 in \cite{BuvoliETDZenodo} along with code for reproducing our numerical experiments. 

In all our numerical experiments, we apply a fine spectral spatial discretization so that error is primarily due to the time integrator. In our first three experiments we impose periodic boundary conditions and solve the PDEs in Fourier space. This is convenient since it leads to an evolution equation of the form (\ref{eq:model_ode_semilinear}) where the matrix $\Lambda$ is diagonal. In our final experiment we consider a more challenging example where $\Lambda$ is a dense matrix. We base our first three numerical experiments from \cite{KassamTrefethen05ETDRK4, grooms2011IMEXETDCOMP} so that our results can be compared with those obtained using other IMEX and ETD schemes.

Since we consider methods of varying order, our experiments are based on the number of function evaluations rather than the step size $h$. We compute reference solutions by using four times as many function evaluations as used in the experiment. To avoid biased results, we average the solutions of at least two convergent methods when forming our reference solutions. For each PDE, we present plots of relative error vs. function evaluations, relative error vs. stepsize, and relative error vs. computational time, where the relative error between two solution vectors $\mathbf{x}$ and $\mathbf{y}$ is $\|\mathbf{x} - \mathbf{y} \|_\infty/ \| \mathbf{x}\|_\infty.$ Though we solve equations in Fourier space, we compute relative errors in physical space. We do not count the time required to initialize ETD coefficients in our time plots. We also make no specific efforts to optimize our code, thus timing results only serve as an indication and may vary under different implementations. The results presented in this paper have been run on a 3.5 Ghz Intel i7 Processor using our double precision Fortran implementation. We describe each of the four problems below.

\begin{figure}[!htb]
\centering

\par{{\bf \noindent Performance Results for Kuramoto-Sivashinsky Equation }}
\vspace{0.5em}
\includegraphics[width=1\linewidth]{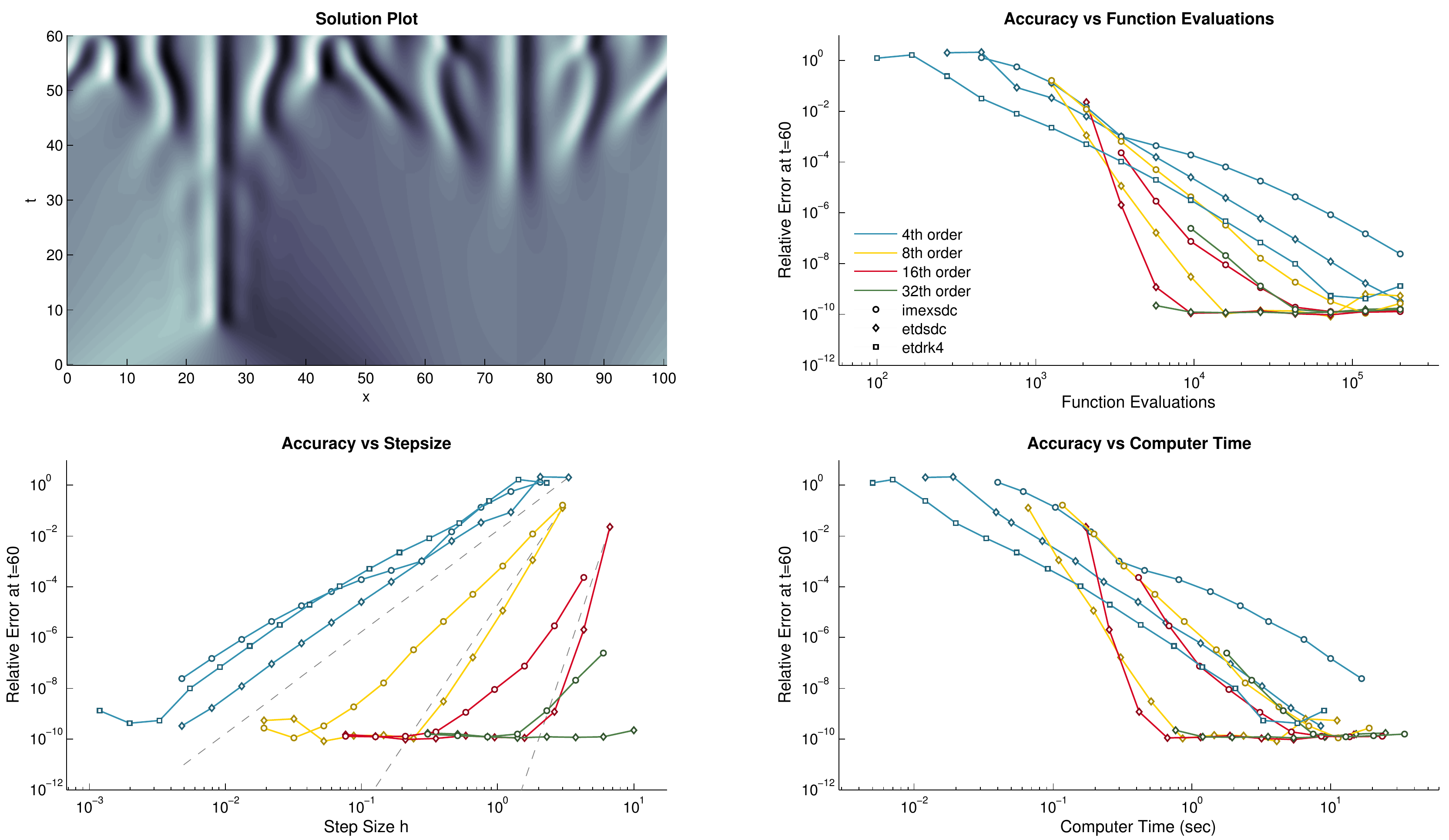}

\vspace{1.5em}
\par{ {\bf \noindent  Performance Results for Nikolaevskiy Equation}}
\vspace{0.5em}
\includegraphics[width=1\linewidth]{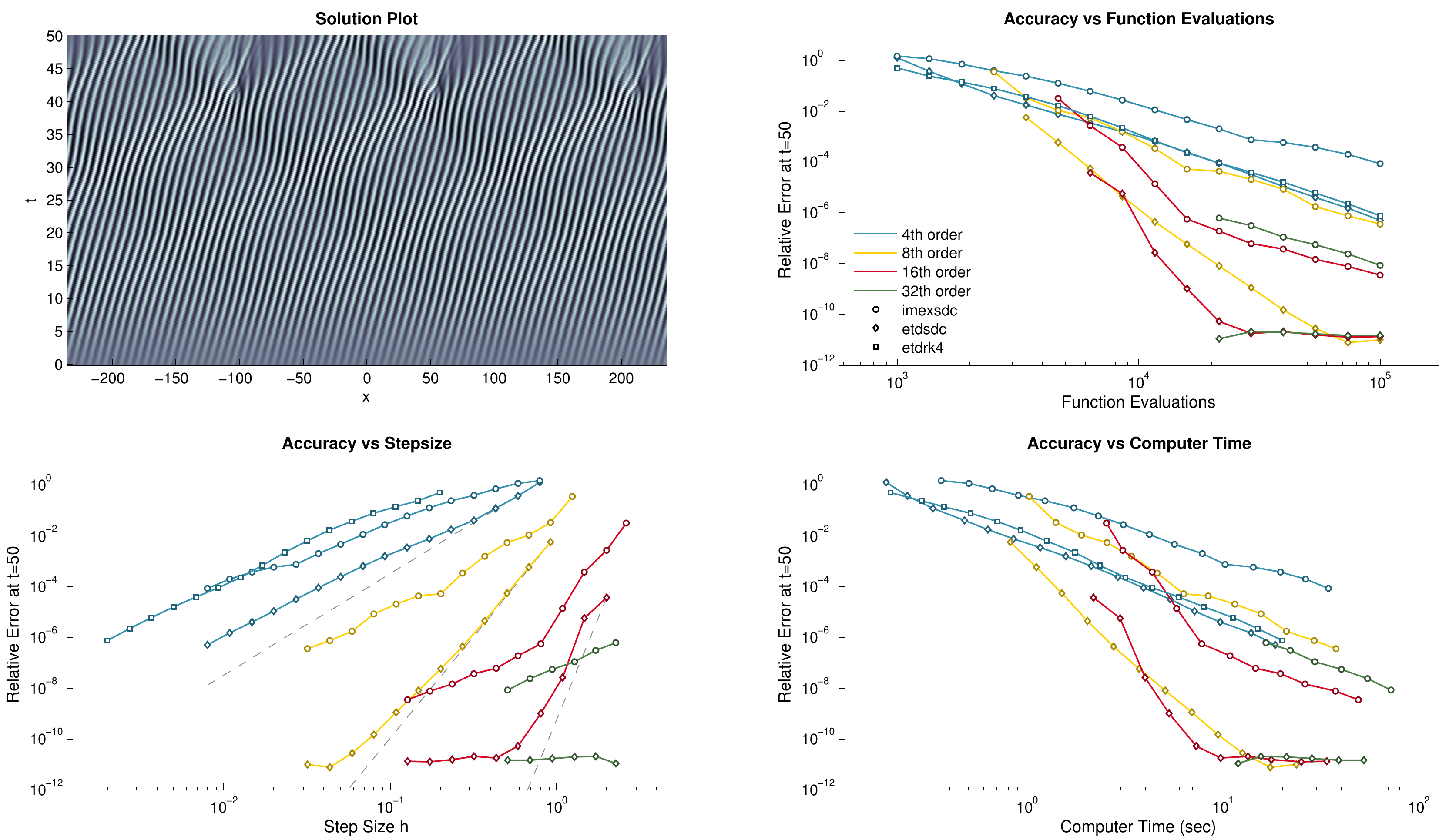}

\caption{Performance results for the Kuramoto-Sivianshi and Nikolaevskiy equations. Gray dashed lines of increasing steepness in the accuracy vs stepsize plots correspond to $O(h^4)$, $O(h^8)$ and $O(h^{16})$, respectively. IMEXSDC schemes experience significant order reduction on both problems.}
\label{fig:results_page1}
\end{figure}

\vskip 1em \noindent
The {\bf Kuramoto-Sivashinsky} (KS) equation models reaction-diffusion systems \cite{kuramoto1976persistent}. As originally presented in \cite{KassamTrefethen05ETDRK4}, we consider the KS equation with periodic boundary conditions:
\begin{align}
& u_t = -u_{xx} -u_{xxxx} - \tfrac{1}{2}\left(u^2\right)_x \label{eq:kuramoto}, \\
& u(x,t=0) = \cos\left( \tfrac{x}{16} \right) \left( 1 + \sin\left( \tfrac{x}{16} \right) \right), \hspace{1em} x\in[0, 64\pi]. \nonumber
\end{align}
We numerically integrate \eref{kuramoto} using a 1024 point Fourier spectral discretization in $x$ and run the simulation out to $t=60$. The KS equation has a dispersive linear term $\Lambda$ with eigenvalues given by $\lambda(k) = k^2 - k^4$, where $k$ denotes the Fourier wavenumber. We present our numerical results in Figure \ref{fig:results_page1}.

\vskip 1em \noindent
The {\bf Nikolaevskiy} equation was originally developed for studying seismic waves \cite{Nikolaevskiy} and now serves as a model for pattern formation in a variety of systems \cite{simbawa2010nikolaevskiy}. As originally presented in \cite{grooms2011IMEXETDCOMP}, we consider the Nikolaevskiy equation with periodic boundary conditions:
\begin{align}
& u_t = \alpha \partial^3_x u + \beta \partial^5_x u -\partial_x^2 \left( r - ( 1 + \partial^2_x)^2 \right) u -\tfrac{1}{2}\left( u^2 \right)_x \label{eq:nikolaevskiy}, \\
& u(x,t=0) = \sin(x) + \epsilon \sin(x/25), \hspace{1em} x\in[-75\pi, 75\pi] \nonumber
\end{align}
where $r=1/4$, $\alpha = 2.1$, $\beta = 0.77$, and $\epsilon = 1/10$.  We solve the Nikolaevskiy equation using a 4096 point Fourier spectral discretization in $x$ and run the simulation out to $t=50$. The Nikolaevskiy equation has a dissipative and dispersive linear term with eigenvalues given by ${\lambda(k) = k^2(r - (1-k^2)^2) - i\alpha k^3 + i\beta k^5}$, where $k$ denotes the Fourier wavenumber. We present our numerical results in Figure \ref{fig:results_page1}.

\begin{figure}[!htb]
\centering

\par{{\bf \noindent Performance Results for Quasigeostrophic Equation}}
\vspace{0.5em}

\includegraphics[width=1\linewidth]{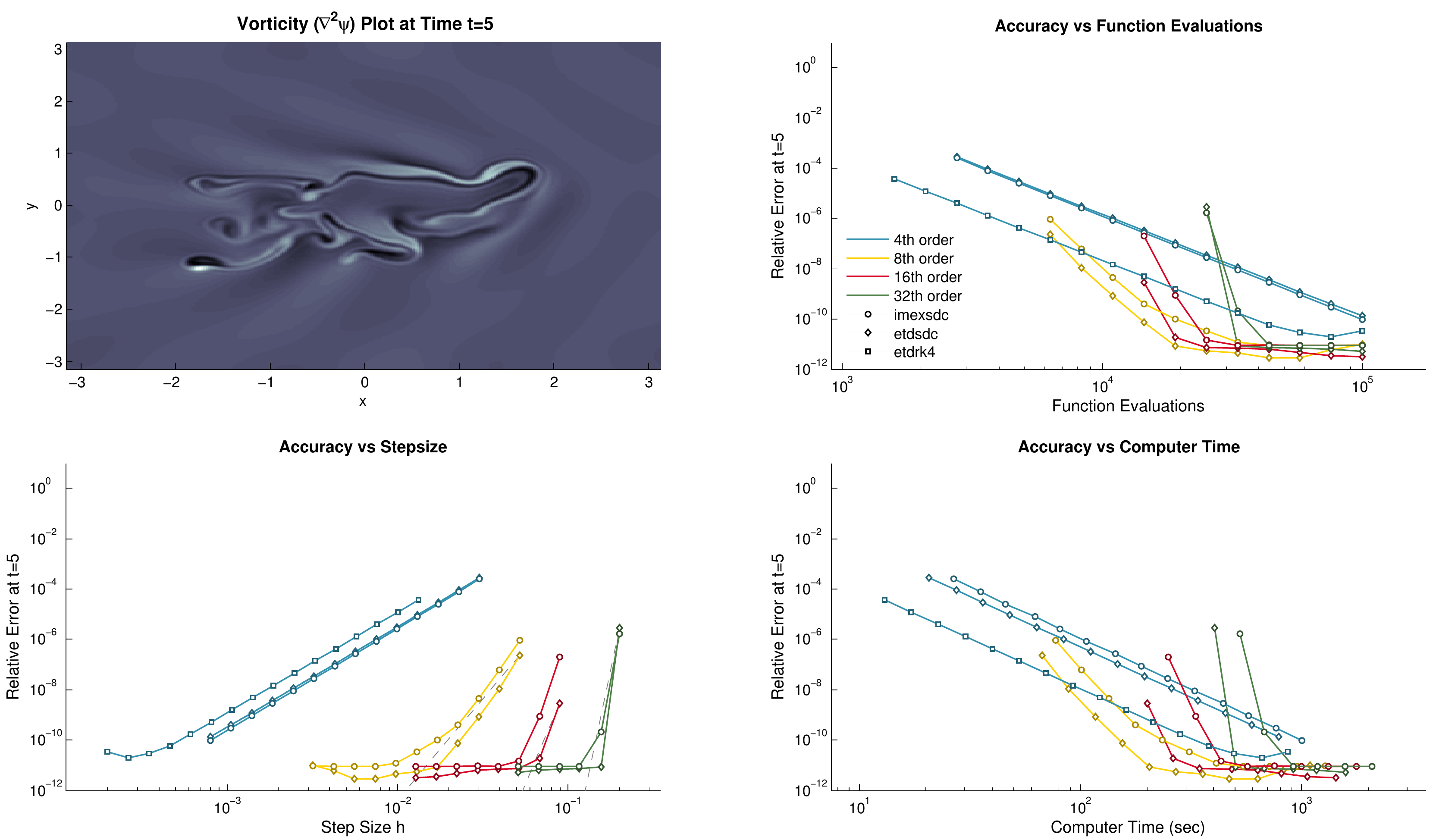}

\vspace{1.5em}
\par{{\bf Performance Results for Korteweg-de Vries Equation}}
\vspace{0.5em}

\includegraphics[width=1\linewidth]{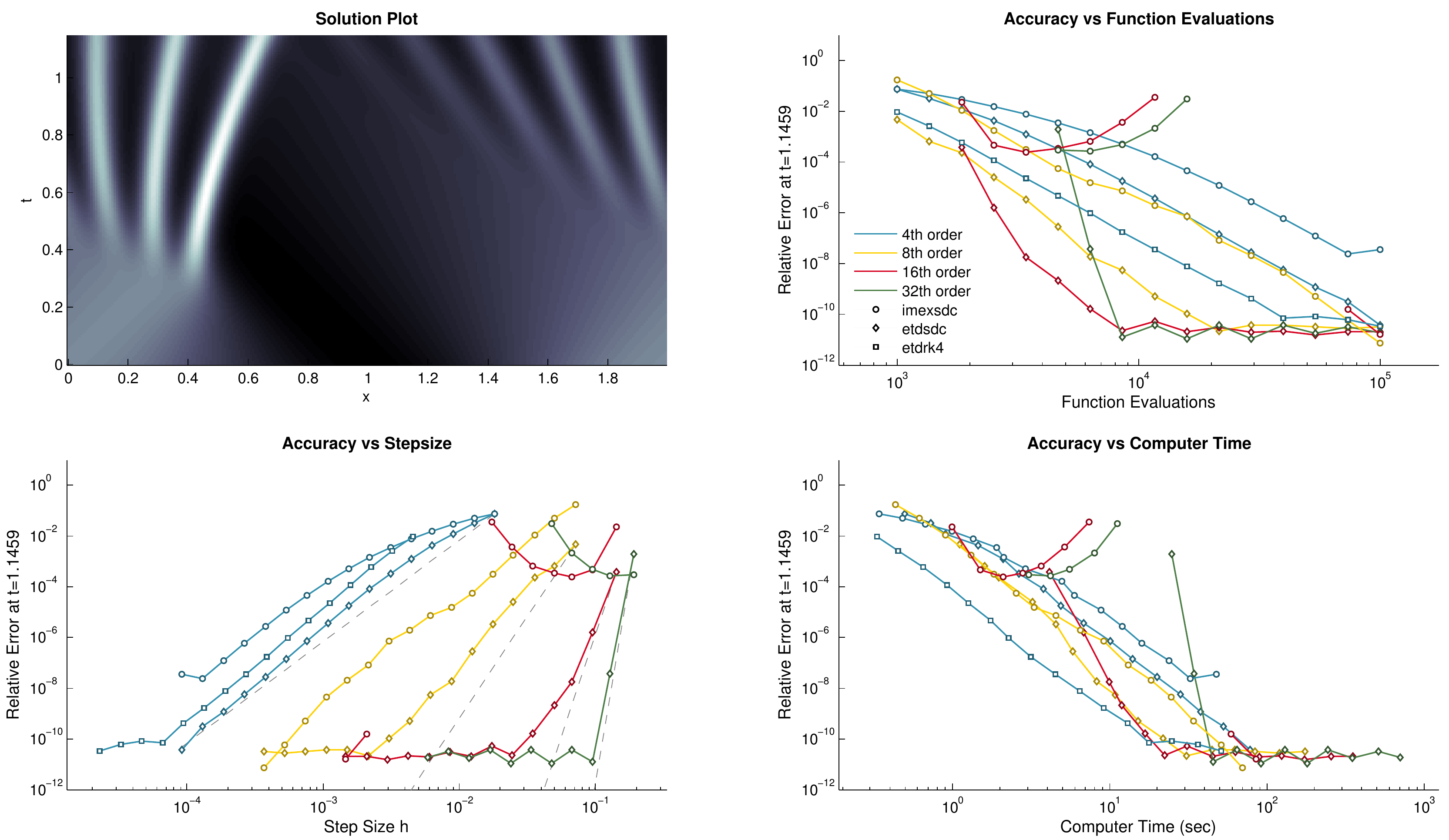}

\caption{Performance results for the Quasigeostrophic and Korteweg-de Vries equations. Dashed lines of increasing steepness in the accuracy vs stepsize plots correspond to $O(h^4)$, $O(h^8)$ $O(h^{16})$ and $O(h^{32})$, respectively. Notice that high-order IMEXSDC schemes are unstable on the KDV equation. Order reduction does not occur for any method on the quasigeostrophic equation, but affects both IMEXSDC and ETDSDC schemes on the KDV equation.}
\label{fig:results_page2}
\end{figure}

\vskip 1em \noindent
The {\bf quasigeostrophic} (QG) equations model a variety of atmospheric and oceanic phenomena \cite{pedloskygeophysical}. As originally presented in \cite{grooms2011IMEXETDCOMP}, we consider the barotropic QG equation on a $\beta$-plane with linear Ekman drag and hyperviscous diffusion of momentum with periodic boundary conditions,
\begin{align}
& \partial_t \nabla^2 \psi  = -\left[\beta \partial_x \psi  + \epsilon \nabla^2 \psi + \nu \nabla^{10} \psi + \mathbf{u} \cdot \nabla (\nabla^2 \psi) \right] \label{eq:quasigeostrophic} \\
& \psi(x,y,t=0) = \frac{1}{8} \exp\left(-8\left(2y^2+x^2/2 - \pi/4 \right)^2 \right), \nonumber \\
& (x,y) \in[-\pi,\pi] \nonumber
\end{align}
where $\psi(x,y)$ is the stream function for two-dimensional velocity ${\mathbf{u} = (-\partial_y \psi, \partial_x \psi)}$, $\epsilon = 1/100$,  and $\nu = 10^{-14}$. We run the simulation to time $t=5$ using a $256\times 256$ point Fourier discritization. We consider a different initial condition than the one presented in \cite{grooms2011IMEXETDCOMP}, since $\nabla^2 \psi(x,y)$  was originally chosen to be discontinuous at the point $(0,0)$. We note that \eref{quasigeostrophic} describes the change in the vorticity $\omega = \nabla^2 \psi$ in terms of the stream function $\psi$. In order to obtain $\psi$ at each timestep, it is necessary to solve Poisson's equation $\nabla^2 \psi = \omega$. Since we are solving in Fourier space, it follows that 
\begin{equation*}
\hat{\psi}_{k,l} = 
\begin{cases}
0 & k=l=0 \\
-\frac{\hat{\omega}}{k^2 + l^2} & \text{otherwise}
\end{cases}
\end{equation*}
where $k$ and $l$ are the Fourier wave numbers and $\hat{\psi}$, $\hat{\omega}$ denote the discrete Fourier transforms of $\psi$ and $\omega$. The QG equation has a linear term with strong dissipation and mild dispersion with eigenvalues given by $\lambda(k,l) = \frac{-ik - \epsilon k^2}{k^2 + l^2} - \nu(k^8 + l^8)$. We present our numerical results in Figure \ref{fig:results_page2}.

\clearpage


\vskip 1em \noindent
The {\bf Korteweg-de Vries} (KDV) equation describes weakly nonlinear shallow water waves. In 1965 Kruskal and Zabusky observed that smooth initial conditions could give rise to soliton solutions \cite{zabusky1965interaction}. As in their original numerical experiment, we consider the KDV equation on a periodic domain
\begin{align*}
& u_t = -\left[ \delta u_{xxx} + \tfrac{1}{2}(u^2)_x \right] \\
& u(x,t=0) = \cos(\pi x), \hspace{1em} x \in [0,2]
\end{align*}
where $\delta = 0.022$ and the simulation is run out to time $t = 3.6/\pi$. The eigenvalues of the linear terms are given by $\lambda(k,l) = \delta ik^3$; thus this equation possess a purely dispersive linear term. Unlike our previous examples, we solve this PDE in physical space where the resulting differentiation matrix is no longer diagonal. The nondiagonal case is more challenging since the coefficients $w_{i,n}$ in \eref{W_lagrange} are now matrix functions. In practice it would be more efficient to consider a lower-order spatial description and apply Krylov space or contour integral techniques that avoid explicitly initializing the requisite ETD matrices. Nevertheless, we consider this example to test the robustness of the scaling and squaring algorithm. For IMEXSDC$^M_N$ schemes it is necessary to repeatedly solve the system $\Lambda x = f$ at each timestep. We perform an initial $LU$ factorization of $\Lambda$ to expedite this process. We present our numerical results for the KDV equation in Figure \ref{fig:results_page2}. 


\subsection{Discussion}

Our results demonstrate that high-order methods can lead to significant speedup when solving nonlinear wave equations to high accuracy. Methods of order 8 and 16 were able to achieve the smallest error using the fewest function evaluations and the least overall CPU time. Interestingly, the error threshold separating good and bad performance for high and low order methods varied significantly in each experiment. Overall, ETDSDC methods consistently achieved better accuracy than corresponding IMEXSDC methods, and did not suffer from crippling order reduction on any of the problems we tested. Amongst the fourth order methods, ETDRK4 is more efficient than either ETDSDC$_4^3$ or IMEXSDC$_4^3$. Moreover, ETDRK4 is the fastest method for computing solutions if error tolerances are large. Methods of order 32 were generally less competitive than those of 8th or 16th order, and should only be considered in situations where extreme precision is necessary and quad/arbitrary-precision arithmetic allow for relative errors significantly below ${1 \times 10^{-12}}$. Finally, for diagonal $\Lambda$, the time required to initialize the ETD coefficients was insignificant as compared to overall computational time even for 32nd order method.


High-order ETDSDC methods continued to perform well even in the non-diagonal case, and we found no evidence of catastrophic roundoff error when forming the ETD matrix coefficients $w_{i,l}(h_i \Lambda)$. For nondiagonal $\Lambda$, high-order ETDSDC$_N^M$ schemes require large amounts of memory and time to initialize and store the $N^2-N$ requisite matrices. Moreover, the expensive matrix multiplications at each timestep reduced their overall competitiveness. To improve the performance of ETDSDC schemes on higher dimensional problems with non-diagonal linear operators, it becomes essential to use techniques that avoid explicitly storing the ETD matrices. 

High-order IMEXSDC schemes were unstable when solving the KDV equation on fine grids in both physical and Fourier space. Through additional numerical testing we find that IMEXSDC schemes can be unstable when integrating other nonlinear wave equations with dispersive linear terms such as the nonlinear Schr\"{o}dinger equation. 

We make several additional comments regarding our numerical experiments. The benefits of using high-order methods is greatly reduced if the initial conditions are not smooth, though in certain situations we found that high-order methods are rendered no less efficient than lower-order counterparts. The size of the integration window also affects the difference in performance between high and low-order methods, with the high-order methods generally benefiting on larger time domains. Chaotic equations can cause additional complications, as small perturbations due to rounding errors grow exponentially and contaminate overall accuracy. This was the case for the the KS equation where we were not able to integrate further without damaging the quality of the reference solution.


\section{Conclusion}

We have demonstrated that high-order ETD spectral deferred correction schemes possess excellent accuracy/stability properties and outperform existing ETD and IMEX methods when solving nonlinear wave equations to high accuracy. Our proposed methodology for initializing ETD coefficients is robust and can be successfully applied to ETDSDC schemes up to 32 order accuracy, even for equations with non-diagonal linear operator $\Lambda$. We have also highlighted the advantages of ETD spectral deferred correction methods as compared with IMEXSDC schemes. Our new ETD schemes consistently outperform their IMEX counterparts, do not appear to suffer from crippling order reduction, and retain stability on equations with dispersive linear terms. 
\includegraphics[height=0.75em]{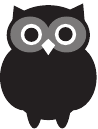}

\section*{Acknowledgements}
I would like to thank Randall J. LeVeque for the many useful discussions over the course of this project and for his comments on drafts of this work. I would also like to thank Kristina Callaghan for her help editing the manuscript. This research was supported in part by funding from the Applied Mathematics Department at the University of Washington and NSF grant DMS-1216732.


\appendix

\section{Choosing $\delta(m)$}
\label{ap:A1}

\noindent We describe a simple choice for $\delta(m)$ from Table \ref{tab:phi_methods}; more sophisticated alternatives are developed in \cite{koikari2007error}. Let $\Lambda$ be a matrix or scalar and let $\mathbf{A} = \Lambda/2^s$. We seek an integer $s$ so that $\varphi_n(\mathbf{A})$ can be initialized via its $m$th order Taylor series without admitting an error larger than $\epsilon$, assuming exact arithmetic. $\varphi_n(\mathbf{A})$ can be approximated by
\begin{equation*}
\varphi_n^{m}(\mathbf{A}) = \sum_{k=0}^m \frac{\mathbf{A}^{k}}{(k+n)!} + \BigO{\|\mathbf{A}\|^{m+1}}.
\end{equation*}
The error $E_n(\mathbf{A}) = \| \varphi_n(\mathbf{A}) - \varphi^m_n(\mathbf{A}) \|$ can be expressed as
\begin{equation*}
E_n(\mathbf{A}) = \left| \left| \sum_{k={m+1}}^\infty \frac{\mathbf{A}^{k}}{(k+n)!} \right| \right| \le  \sum_{k={m+1}}^\infty  \frac{ \| \mathbf{A} \|^k}{(k+n)!}.
\end{equation*}
For any $n$, $E_n(\mathbf{A}) $ can be bounded above by
\begin{equation*}
\sum_{k=1}^\infty \frac{\|\mathbf{A} \|^{k+m}}{k!(m+1)!} = \frac{\|\mathbf{A} \|^{m+1} \exp(\|\mathbf{A} \|)}{(m+1)!}.
\end{equation*}
Assuming exact arithmetic, we can guarantee that $E_n(\mathbf{A}) < \epsilon$ for any $s$ so that $\mathbf{A} = \Lambda/2^s$ satisfies
\begin{equation*}
\epsilon  = \frac{\|\mathbf{A} \|^{m+1} \exp(\|\mathbf{A} \|)}{(m+1)!}.
\end{equation*}
To avoid solving a nonlinear system for each $m$ and $\epsilon$, we fix $m=20$, $\epsilon=1\times 10^{-16}$ and let $\|\mathbf{A}\|\le \rho$. Solving the cooresponding nonlinear equation for $\rho$, leads to $\rho = 1.4 \approx 1.0$. Therefore our condition on $\mathbf{A}$ reduces to $\| \mathbf{A} \| = \| \Lambda/2^s \| \le 1$.
In our numerical codes we choose $\| \mathbf{A} \| = \max\{\|\mathbf{A}\|_1,\|\mathbf{A}\|_{\infty}\}$. This leads to the condition
\begin{equation*}
s = \max \left \{ 0, \frac{\ln \left( \max\{\|\mathbf{A}\|_1,\|\mathbf{A}\|_{\infty}\} \right)}{\ln(2)} \right \}.
\end{equation*}


\bibliographystyle{plain}
\bibliography{references_etd,references_sdc,references_other}

\end{document}